\documentclass[aap,MSNbibl,citesort,seceqn,dvips]{arximspdf}
\usepackage{graphicx}
\doi{10.1214/11-AAP808} 
\volume{22}
\issue{4}
\pubyear{2012}
\firstpage{1650}
\lastpage{1692}

\makeatletter

\newtheorem{thmm}{Theorem}[section]
\newtheorem{prop}[thmm]{Proposition}
\newtheorem{lem}[thmm]{Lemma}
\newproclaim{defn}[thmm]{Definition}
\newproclaim{rem}[thmm]{Remark}
\def\N{\mathbb{N}}
\def\Z{\mathbb{Z}}
\def\R{\mathbb{R}}
\def\P{\mathbb{P}}
\def\E{\mathbb{E}}
\newcommand{\bzd}{{(\Z^d)^*}}
\newcommand{\bzdl}{\Lambda^*}
\newcommand{\C}{\mathcal{C}}
\newcommand{\A}{\mathcal{A}}
\newcommand{\ormd}{\mathrm{d}}
\newcommand{\rmd}{\mathrm{d}}
\newcommand{\bvar}{\operatorname{\mathbb{V}ar}}

\def\sb{{\subset}}

\newcommand{\zd}{{\mathbb{Z}^d}}
\newcommand{\RR}{\mathbb{R}}

\def\b{\beta}

\def\d{\delta}

\def\phi{\varphi}

\def\g{\gamma}

\def\s{\sigma}

\def\t{\tau}

\def\x{\xi}

\def\L{\Lambda}

\makeatother

\begin{document}
\begin{frontmatter}

\title{Existence of random gradient states\thanksref{T1}}
\runtitle{Existence of random gradient states}
\thankstext{T1}{Support by the TUM Institute
for Advanced Study (TUM-IAS), by the Sonderforschungsbereich
SFB | TR12---Symmetries and Universality in Mesoscopic Systems,
and by the University of Bochum.}

\begin{aug}
\author[A]{\fnms{Codina} \snm{Cotar}\corref{}\ead[label=e1]{ccotar@fields.utoronto.ca}}
\and
\author[B]{\fnms{Christof} \snm{K\"ulske}\ead[label=e2]{Christof.Kuelske@ruhr-uni-bochum.de}}
\runauthor{C. Cotar and C. K\"ulske}
\affiliation{The Fields Institute and Ruhr-University of Bochum}
\address[A]{The Fields Institute\\
222 College Street\\
Toronto, Ontario\\
Canada M5T 3J1\\
\printead{e1}} 
\address[B]{Fakult\"at f\"ur Mathematik\\
Ruhr-University of Bochum\\
Postfach 102148\\
44721, Bochum\\
Germany\\
\printead{e2}}
\end{aug}

\received{\smonth{8} \syear{2010}}
\revised{\smonth{9} \syear{2011}}

%
\begin{abstract}
We consider two versions of random gradient models.
In model A the interface feels a bulk term of random fields
while in model B the disorder enters through the potential acting on
the gradients. It is well known that for gradient models without disorder
there are no Gibbs measures in infinite-volume in dimension
$d = 2$, while there are ``gradient Gibbs measures'' describing an
infinite-volume distribution for the gradients of the field, as was
shown by Funaki and Spohn. Van Enter and K\"ulske proved that adding a
disorder term as in model A prohibits
the existence of such gradient Gibbs measures for general interaction
potentials in $d = 2$.

In the present paper we prove the existence of shift-covariant gradient
Gibbs measures with a given tilt $u\in\R^d$ for model A when $d\geq3$
and the disorder has mean zero, and for model B when $d\geq1$. When
the disorder has nonzero mean in model A, there are no shift-covariant
gradient Gibbs measures for $d\ge3$. We also prove similar results of
existence/nonexistence of the surface tension for the two models and
give the characteristic properties of the respective surface tensions.

\end{abstract}

%

\end{frontmatter}

\section{Introduction} \label{sect:intro}

\subsection{The setup}
Phase separation in $\RR^{d+1}$ can be described by effective interface
models for the study of phase boundaries at a mesoscopic level
in statistical mechanics. Interfaces are sharp boundaries which
separate the different regions of space occupied by different phases.
In this class of models, the interface is modeled as the graph of a random
function from $\Z^d$ to $\Z$ or to $\RR$ (discrete or continuous
effective interface models). For background and earlier results on
continuous and discrete interface models without disorder, see, for
example,~\cite{BY,CD,CDM,DGI,FL,FS,gos} and references therein. In our setting, we will
consider the case of continuous interfaces with disorder as introduced
and studied previously in~\cite{EK} and~\cite{KO}. Note also that
discrete interface models in the presence of
disorder have been studied, for example, in~\cite{BK1} and~\cite{BK2}.
We will introduce next our two models of interest.

In our setting, the fields $\phi(x)\in\RR$ represent height variables
of a random interface at the site $x\in\zd$. Let $\Lambda$ be a finite
set in $\Z^d$ with boundary
%
%
\begin{eqnarray}
\partial\Lambda:=\{x\notin\Lambda, \Vert x-y\Vert =1 \mbox{ for some } y\in
\Lambda
\}
\nonumber
\\[-8pt]
\\[-8pt]
\eqntext{\mbox{where } \|x-y\|=\displaystyle\sum^d_{i=1}|x_i-y_i|.}
\end{eqnarray}
On the boundary we set a boundary condition $\psi$ such that $\phi
(x)=\psi(x)$ for $x\in\partial\Lambda$. Let $(\Omega, \mathcal{F},\P
)$ be
a probability
space; this is the probability space of the disorder, which will be
introduced below. We denote by the symbol $\E$ the expectation w.r.t.~$\P$.

Our two models are given in terms of the \textit{finite-volume
Hamiltonian} on~$\Lambda$.
\begin{longlist}[(A)]
\item[(A)] For model A the Hamiltonian is
%
%
\begin{eqnarray}\qquad
\label{eqn00}
H_{\Lambda}^\psi[\xi](\phi)&:=&\frac{1}{2}\mathop{\sum_{x, y\in\Lambda
}}_{|x-y|=1}V\bigl(\phi(x)-\phi(y)\bigr)+\mathop{\sum_{x\in\Lambda,y\in\partial\Lambda
}}_{
|x-y|=1}V\bigl(\phi(x)-\psi(y)\bigr)
\nonumber
\\[-8pt]
\\[-8pt]
\nonumber
&&{}+ \frac{1}{2}\sum_{x\in\Lambda}\xi
(x)\phi(x),
\end{eqnarray}
where the random fields $(\xi(x))_{x\in\Z^d}$ are assumed to be i.i.d.
real-valued random variables, with finite nonzero second moments. The
disorder configuration $(\xi(x))_{x\in\Z^d}$ denotes an
arbitrary fixed configuration of external fields, modeling a
``quenched'' (or
frozen) random environment. We assume that $V\in C^2(\Bbb R)$ is an
even function with quadratic growth at infinity:
%
%
\begin{equation}
\label{tag2}
V(s)\ge As^2-B,\qquad s\in\Bbb R,
\end{equation}
for some $A>0, B\in\Bbb R$. We assume also that there exists $C_2>0$
such that
%
%
\begin{equation}
\label{tag22}
V''(s)\le C_2 \qquad \mbox{for all }  s\in\RR.
\end{equation}
\item[(B)] For each bond $(x,y)\in\Z^d\times\Z^d,|x -y|=1$, we define
the measurable map $V_{(x,y)}^{\omega}(s)\dvtx (\omega,s)\in\Omega\times
\RR
\rightarrow\RR$. Then $V_{(x,y)}^\omega$ is a random real-valued
function and $V_{(x,y)}^\omega$ are assumed to be i.i.d. random
variables as $(x,y)$ ranges over the bonds.
Let $B_{(x,y)}^\omega$ be a family of i.i.d. real-valued random
variables with $\E|B_{(x,y)}|<\infty$.

We assume that for some given $A,C_2>0$, $V_{(x,y)}^\omega$ obey for
$\P
$-almost every $\omega\in\Omega$ the following bounds, uniformly in the
bonds $(x,y)$:
%
%
\begin{equation}
\label{tag23}
A s^2 -B_{(x,y) }^\omega\leq V_{(x,y)}^\omega(s)\leq C_2 s^2 \qquad\mbox{for
all } s\in\RR.
\end{equation}
We assume also that for each fixed $\omega\in\Omega$ and for each bond
$(x,y)$, $V_{(x,y)}^\omega\in C^2(\Bbb R)$ is an even function. Then
for model B we define the Hamiltonian for each fixed $\omega\in\Omega
$ by
%
%
\begin{eqnarray}
\label{eqn000}
H_{\Lambda}^\psi[\omega](\phi)&:=&\frac{1}{2}\sum_{x,y\in\Lambda
,|x-y|=1}V_{(x,y)}^\omega\bigl(\phi(x)-\phi(y)\bigr)
\nonumber
\\[-9pt]
\\[-9pt]
\nonumber
&&{}+\sum_{x\in\Lambda,y\in
\partial\Lambda,|x-y|=1}V_{(x,y)}^\omega\bigl(\phi(x)-\psi(y)\bigr).
\end{eqnarray}
\end{longlist}

The two models above
are prototypical ways to add randomness which preserve the gradient
structure, that is, the Hamiltonian depends only on the gradient field
$(\phi(x)-\phi(y))_{x,y\in\Z^d, |x-y|=1}$. Note that for $d=1$
our interfaces can be used to model a polymer chain; see, for example,
\cite{denholl}. Disorder in the Hamiltonians models impurities in the
physical system. Models A and B can be regarded as modeling two
different types of impurities, one affecting the interface height, the
other affecting the interface gradient.

The rest of the \hyperref[sect:intro]{Introduction} is structured as follows: in Section
\ref{AG} we define in detail the notions of finite- and infinite-volume
(gradient) Gibbs measures for model A, in Section~\ref{BG} we sketch
the corresponding notions for model B, in Section~\ref{ST} we
introduce the notion of surface tension for the two models, and in
Section~\ref{MR} we present our main results and their connection to
the existing literature.\vspace*{-3pt}

\subsection{Gibbs measures and gradient Gibbs measures for model A}\label{AG}\vspace*{-3pt}

\subsubsection{\texorpdfstring{$\phi$-Gibbs measures}{phi-Gibbs measures}}

Let $C_b(\R^{\zd})$ denote the set of continuous and bounded functions
on $\R^{\zd}$. The functions considered are functions of the interface
configuration $\phi$, and
continuity is with respect to each coordinate $\phi(x),x\in\Z^d,$ of
the interface. For a finite region $\Lambda\subset\Z^d$, let $\rmd
\phi
_{\Lambda}:=\prod_{x\in\Lambda}\,\rmd\phi(x)$ be the Lebesgue measure
over $\R^{\Lambda}$.

Let us
first consider model A only, and let us define the $\phi$-Gibbs
measures for \textit{fixed} disorder $\xi$.\vspace*{-3pt}
\begin{defn}[(Finite-volume $\phi$-Gibbs measure)]\label{gibbs0}
For a finite region $\Lambda\subset\Z^d$, \textit{the finite-volume Gibbs
measure $\nu_{\Lambda,\psi}[\xi]$ on $\R^{\Z^d}$} with given
Hamiltonian $H[\xi]:=(H_{\Lambda}^\psi[\xi])_{\Lambda\subset\zd,
\psi\in\RR^{\zd}}$, with boundary condition $\psi$ for the field of
height
variables $(\phi(x))_{x\in\Z^d}$ over $\Lambda$, and with a fixed
disorder configuration~$\xi$, is defined by
%
%
\begin{equation}
\label{tag0'}
\nu_{\Lambda}^\psi[\xi](\ormd\phi):=\frac{1}{Z_{\Lambda}^\psi
[\xi]}\exp
\{-H_{\Lambda}^\psi[\xi](\phi)\}\,\rmd\phi_\Lambda\delta_\psi
(\ormd\phi_{{\Z}^d\setminus\Lambda}),\vadjust{\goodbreak}
\end{equation}
where
\[
Z_{\Lambda}^\psi[\xi]:=\int_{{\R}^{\Z^d}}\exp\{-H_{\Lambda
}^\psi
[\xi](\phi)\}\,\rmd\phi_\Lambda\delta_\psi(\ormd\phi_{\Z
^d\setminus
\Lambda})
\]
and
\[
\delta_\psi(\ormd\phi_{\Z^d\setminus\Lambda}):=\prod_{x\in\Z
^d\setminus
\Lambda}\delta_{\psi(x)}(\ormd\phi(x)).
\]
\end{defn}

It is easy to see that the growth
condition on $V$ guarantees the finiteness of the integrals appearing
in (\ref{tag0'}) for
all arbitrarily fixed choices of $\xi$.
\begin{defn}[($\phi$-Gibbs measure on $\zd$)]
\label{gibbs}
The probability measure $\nu[\xi]$ on $\RR^{\zd}$ is called an \textit{(infinite-volume) Gibbs measure} for the $\phi$-field with given
Hamiltonian $H[\xi]:=(H_{\Lambda}^\psi[\xi])_{\Lambda\subset\zd,
\psi
\in
\RR^{\zd}}$ ($\phi$-Gibbs measure for short), if it satisfies the
DLR equation
%
%
\begin{equation}
\label{dlrgibbs}
\int\nu[\xi](\ormd\psi)\int\nu_{\Lambda}^\psi[\xi](\ormd\phi
)F(\phi
)=\int\nu[\xi](\ormd\phi)F(\phi),
\end{equation}
for every finite $\Lambda\subset\zd$ and for all $F\in C_b(\R^{\zd})$.
\end{defn}

We discuss next the case of interface models without disorder, that is,
with $\xi(x)=0$ for all $x\in\Z^d$ in model A.
Let $\nu^\psi_\Lambda[\xi=0],\L\in\Z^d$, denote the
finite-volume Gibbs
measure for $\L$ and with boundary condition $\psi$. Then an
infinite-volume Gibbs measure $\nu[\xi=0]$ exists under condition
(\ref
{tag2}) only when $d\ge3$, but not for $d=1,2$, where the field
``delocalizes'' as $\Lambda\nearrow\Z^d$ (see~\cite{FP}).

In the case of interfaces with disorder as in model A, it has been
proved in~\cite{KO} that the $\phi$-Gibbs measures do not exist when
$d=2$. A similar argument as in~\cite{KO} can be used to show that
$\phi
$-Gibbs measures do not exist for model A when $d=1$.

\subsubsection{\texorpdfstring{$\nabla\phi$-Gibbs measures}{nabla phi-Gibbs measures}}

We note that the Hamiltonian $H_\L^\psi[\xi]$ in mod\-el~A, respectively,
$H_\L^\psi[\omega]$ in model B, changes only by a configuration-independent
constant under the joint shift $\phi(x)\rightarrow\phi(x)+c$ of all
height variables $\phi(x),x\in\Z^d,$
with the same $c\in\RR$. This holds true for any fixed configuration
$\xi$, respectively,~$\omega$. Hence, finite-volume
Gibbs measures transform under a shift of the boundary condition
by a shift of the integration variables. Using this invariance under height
shifts, we can lift the finite-volume measures to measures on gradient
configurations, that is, configurations of height differences across
bonds,
defining the gradient finite-volume Gibbs measures. Gradient Gibbs
measures have the advantage that they may exist, even
in situations where the Gibbs measure does not. Note that the concept
of $\nabla\phi$-measures is general and does not refer only to the
disordered models. For example, in the case of interfaces without
disorder $\nabla\phi$-Gibbs measures exist for all $d\ge1$.

We next introduce the \textit{bond variables on $\zd$}. Let
\[
\bzd:=\{b=(x_b,y_b) | x_b,y_b\in\zd,\|x_b-y_b\|=1,b \mbox{ directed
from } x_b \mbox{ to } y_b\};
\]
note that each undirected bond appears twice in $\bzd$.
For $\phi=(\phi(x))_{x\in\zd}$ and $b=(x_b,y_b)\in\bzd$, we
define the
\textit{height differences} $\nabla\phi(b):=\phi(y_b)-\phi(x_b)$.
The height variables $\phi=\{\phi(x)\dvtx x\in\zd\}$ on $\zd$ automatically
determine a field of height differences
$\nabla\phi=\{\nabla\phi(b)\dvtx b\in\bzd\}$. One can therefore
consider the
distribution $\mu$ of $\nabla\phi$-field
under the $\phi$-Gibbs measure $\nu$. We shall call $\mu$ the
$\nabla
\phi$-Gibbs measure. In fact, it is possible to define
the $\nabla\phi$-Gibbs measures directly by means of the DLR equations
and, in this sense, $\nabla\phi$-Gibbs measures exist for
all dimensions $d\ge1$.

A sequence of bonds $\C=\{b^{(1)},b^{(2)},\ldots,b^{(n)}\}$ is called a
\textit{chain} connecting~$x$ and $y$, $x,y\in\zd$, if
$x_{b_1}=x,y_{b^{(i)}}=x_{b^{(i+1)}}$ for $1\le i\le n-1$ and
$y_{b^{(n)}}=y$. The chain is called a \textit{closed loop} if
$y_{b^{(n)}}=x_{b^{(1)}}$. A \textit{plaquette} is a closed loop $\A
=\{
b^{(1)},b^{(2)},b^{(3)},b^{(4)}\}$ such that $\{x_{b^{(i)}},i=1,\ldots
,4\}$ consists of
four different points.

The field $\eta=\{\eta(b)\}\in\RR^{\bzd}$ is said to satisfy
\textit
{the plaquette condition} if
%
%
\begin{eqnarray}
\eta(b)&=&-\eta(-b) \qquad\mbox{for all } b\in\bzd\quad \mbox{and}
\nonumber
\\[-8pt]
\\[-8pt]
\nonumber
\sum_{b\in
\A}\eta
(b)&=&0\qquad \mbox{for all plaquettes } \A \mbox{ in } \zd,
\end{eqnarray}
where $-b$ denotes the reversed bond of $b$. Let
%
%
\begin{equation}
\chi=\bigl\{\eta\in\RR^{(\zd)^*} \mbox{which satisfy the plaquette
condition}\bigr\}
\end{equation}
and let $L_r^2, r>0$, be the set of all $\eta\in\RR^{\bzd}$ such that
\[
|\eta|^2_r:=\sum_{b\in\bzd}|\eta(b)|^2e^{-2r\|x_b\|}<\infty.
\]
We denote $\chi_r=\chi\cap L_r^2$ equipped with the norm $|\cdot|_r$.
For $\phi=(\phi(x))_{x\in\zd}$ and $b\in\bzd$, we define $\eta
(b):=\nabla\phi(b)$. Then $\nabla\phi=\{\nabla\phi(b)\dvtx b\in\bzd\}$
satisfies the plaquette condition. Conversely, the heights $\phi^{\eta
,\phi(0)}\in\RR^\zd$ can be constructed from height differences
$\eta$
and the
height variable $\phi(0)$ at $x=0$ as
%
%
\begin{equation}
\label{19}
\phi^{\eta,\phi(0)}(x):=\sum_{b\in\C_{0,x}}\eta(b)+\phi(0),
\end{equation}
where $\C_{0,x}$ is an arbitrary chain connecting $0$ and $x$. Note
that $\phi^{\eta,\phi(0)}$ is well defined if
$\eta=\{\eta(b)\}\in\chi$.

Let $C_b(\chi)$ be the set of continuous and bounded functions on
$\chi
$, where the continuity is with respect to each bond variable $\eta
(b),b\in\bzd$.
\begin{defn}[(Finite-volume $\nabla\phi$-Gibbs measure)]
\label{finvolgrad}
The \textit{finite-volume $\nabla\phi$-Gibbs measure} in $\Lambda$ (or
more precisely, in $\bzdl$) with given Hamiltonian $H[\xi
]:=(H_{\Lambda
}^ \rho[\xi])_{\Lambda\subset\zd, \rho\in
\chi}$, with boundary condition $\rho\in\chi$ and with fixed disorder
configuration $\xi$, is a probability measure $\mu_{\Lambda}^\rho
[\xi]$
on $\chi$ such that for all $F\in C_b(\chi)$, we have
%
%
\begin{equation}
\label{fingradgibbs}
\int_{\chi}\mu_{\Lambda}^\rho[\xi](\ormd\eta)F(\eta)=\int_{\R
^{\zd}}\nu
_{\Lambda}^\psi[\xi](\ormd\phi)F(\nabla\phi),
\end{equation}
where $\psi$ is any field configuration whose gradient field is $\rho$.
\end{defn}

%
\begin{defn}[{[$\nabla\phi$-Gibbs measure on $(\zd)^*$]}]
\label{nablaphigib}
The probability measure~$\mu[\xi]$ on $\chi$ is called an \textit{(infinite-volume) gradient Gibbs measure} with given Hamiltonian $H[\xi
]:=(H_{\Lambda}^\rho[\xi])_{\Lambda\subset\zd, \rho\in\chi}$
($\nabla
\phi$-Gibbs measure for short),
if it satisfies the DLR equation
%
%
\begin{equation}
\label{dlrgrad}
\int\mu[\xi](\ormd\rho)\int\mu_{\Lambda}^\rho[\xi](\ormd\eta
)F(\eta
)=\int\mu[\xi](\ormd\eta)F(\eta),
\end{equation}
for every finite $\Lambda\subset\zd$ and for all $F\in C_b(\chi)$.
\end{defn}

\begin{rem}
Throughout the rest of the paper, we will use the notation $\phi,\psi$
to denote height variables and $\eta,\rho$ to denote gradient variables.
\end{rem}

For $v\in\zd$, we define the shift operators: $\tau_{v}$ for the
heights by $(\tau_{v}\phi)(y):=\phi(y-v) \mbox{ for } y\in\zd \mbox{ and } \phi\in\RR^{\zd}$, $\tau_{v}$ for the bonds by $(\tau
_{v}\eta
)(b):=\eta(b-v)$ for $b\in\bzd \mbox{ and } \eta\in\chi$, and
$\tau_v$
for the disorder configuration by $(\tau_v\xi)(y):=\xi(y-v)$ for
$y\in
\zd$ and $\xi\in\RR^{\zd}$.

 We are now ready to define the main object of interest of
this paper: the \textit{random} (gradient) Gibbs measures.
\begin{defn}[{[Translation-covariant random (gradient) Gibbs
measures for model A]}]
\label{shiftcov1}
A measurable map $\xi\rightarrow\nu[\xi]$ is called \textit{a
translation-covariant random Gibbs
measure} if $\nu[\xi]$ is a $\phi$-Gibbs measure for $\P$-almost every
$\xi$, and if
\[
\int\nu[\tau_v\xi](\ormd\phi)F(\phi)=\int\nu[\xi](\ormd\phi
)F(\tau_v\phi),
\]
for all $v\in\zd$ and for all $F\in C_b(\R^{\zd})$.

A measurable map $\xi\rightarrow\mu[\xi]$ is called \textit{a
translation-covariant random gradient Gibbs measure} if $\mu[\xi]$ is a
$\nabla\phi$-Gibbs measure for $\P$-almost every $\xi$, and if
\[
\int\mu[\tau_v\xi](\ormd\eta)F(\eta)=\int\mu[\xi](\ormd\eta
)F(\tau_v\eta),
\]
for all $v\in\zd$ and for all $F\in C_b(\chi)$.\vadjust{\goodbreak}
\end{defn}

The above notion generalizes the
notion of a translation-invariant (gradient) Gibbs measure to the setup
of disordered
systems.\vspace*{-3pt}

\subsection{Gibbs measures and gradient Gibbs measures for model B}
\label{BG}

The notions of finite-volume (gradient) Gibbs measure and
infinite-volume (gradient) Gibbs measure for model B can be defined
similarly as for model A, with $(V^\omega_{(x,y)})_{(x,y)\in\zd
\times\zd
},\omega\in\Omega$, playing a similar role to $\xi\in\RR^\zd$,
and with~$\omega$ replacing $\xi$ in Definitions~\ref{gibbs0}--\ref
{nablaphigib}. Once we specify the action of the shift
map~$\tau_v$ in this case, we can also define the notion of
translation-covariant random (gradient) Gibbs measure, with $\omega\in
\Omega$ replacing $\xi\in\RR^\zd$ in Definition~\ref{shiftcov1}.

Let $\tau_v,v\in\Z^d,$ be a shift-operator and let $\omega\in
\Omega$ be
fixed. We
will denote by $\nu[\tau_v\omega]$ the infinite-volume Gibbs measure
with given Hamiltonian $\bar{H}[\omega](\phi):=(H_\L^\psi[\omega
](\tau_v\phi))_{\L\subset\Z^d,\psi\in\RR^{\Z^d}}$. This means
that we shift the field of disorded potentials on bonds from
$V_{(x,y)}^\omega$ to $V_{(x+v,y+v)}^\omega$. Similarly, we
will denote by $\mu[\tau_v\omega]$ the infinite-volume gradient Gibbs
measure with given Hamiltonian $\bar{H}[\omega](\eta):=(H_\L^\rho
[\omega](\tau_v\eta))_{\L\subset\Z^d,\rho\in\RR^{\bzd}}$.\vspace*{-3pt}



\subsection{Surface tension}
\label{ST}
The surface tension physically measures the macroscopic
energy of a surface with tilt $u\in\R^d$, that is, a $d$-dimensional
hyperplane located in $\RR^{d+1}$ with normal vector $(-u,1)\in\RR
^{d+1}$. In other words, it measures the free-energy
cost in creating an interface with a given tilt.

Formally, let $\Lambda_N=[-N,N]^d\cap\Z^d,N\in\N$, be a hypercube of
side length \mbox{$2N+1$} with boundary $\partial\Lambda_N.$ We enforce a
fixed tilt $u\in\R^d$ by imposing the boundary condition ${\psi
}_u(x)=x\cdot u$ for $x\in\partial\Lambda_N$. The \textit{finite-volume
surface tension $\sigma_{\Lambda_N}[\xi]$} for model A is then defined
for fixed disorder $\xi$ as
%
%
\begin{eqnarray}
\label{surftensdef}
\sigma_{\Lambda_N}[\xi](u)&:=&-\frac{1}{|\Lambda_N|}\log\int_{\R^\Lambda} \exp(-H_\Lambda
^{\psi
_u}[\xi])\,\ormd\phi_{\Lambda_N}
\nonumber
\\[-8pt]
\\[-8pt]
\nonumber
&=&-\frac{1}{|\Lambda_N|}\log
Z^{{{{\psi
}}_u}}_{\L}[\xi],
\end{eqnarray}
where we recall that $\ormd\phi_{\Lambda_N}:=\prod_{x\in\Lambda
_N}\phi
(x)$. We are interested in the existence and $\xi$-independence of the limit:
\[
\sigma[\xi](u):=\lim_{N\rightarrow\infty}\sigma_{\Lambda_N}[\xi](u).
\]
When it exists, the limit $\sigma[\xi](u)$ is called \textit{(infinite-volume) surface tension}.

For model B the surface tension $\sigma_{\Lambda_N}[\omega](u)$,
respectively, $\sigma[\omega](u)$, is defined similarly, with $\omega
\in
\Omega$ in place of $\xi\in\RR^\zd$, in the above definitions for
model A.\vspace*{-3pt}

\subsection{Main results}
\label{MR}
A main question in interface models is
whether the fluctuations of an interface, that is, restricted to a finite-volume
will remain bounded\vadjust{\goodbreak} when the volume tends to infinity, so that there is
an infinite-volume Gibbs measure (or gradient Gibbs measure) describing a
localized interface. This question is well understood in shift-invariant
continuous interface models without disorder, and it is the purpose of
this paper to study the same question for interface
models with disorder.

When there is no disorder, it is known that the Gibbs measure $\nu[\xi
=0]$ does not exist in infinite-volume for $d=1,2$,
but the gradient Gibbs measure $\mu[\xi=0]$ does exist in
infinite-volume for $d\ge1$. The latter fact is equivalent to saying that
the infinite-volume measure exists constrained on \mbox{$\phi(0)=0$}. On the
question of uniqueness of gradient Gibbs measures, Funaki and Spohn~\cite{FS} showed
that a gradient Gibbs measure is uniquely determined by
the tilt. This result has been extended to a certain class of nonconvex
potentials by Cotar and Deuschel in~\cite{CD}.

For (very) nonconvex $V$, new phenomena appear: There is
a first-order phase transition from uniqueness to nonuniqueness of the
Gibbs measures (at tilt zero), as shown in~\cite{BK} and~\cite{CD}. The
transition is due to the temperature which changes the structure of the
interface. This
phenomenon is related to the phase transition seen in rotator models with
very nonlinear potentials exhibited in~\cite{ESHL1} and~\cite{ESHL2},
where the basic mechanism is an
energy-entropy transition.

What happens in the random models A and B? In~\cite{KO} the authors showed
that for model A there is no disordered
infinite-volume random Gibbs measure for $d = 1,2$. This statement is
not surprising
since there exists no $\phi$-Gibbs measure without
disorder. More surprising is the fact that, as proved in~\cite{EK}, for
model A there is also no disordered shift-covariant gradient Gibbs
measure when $d=1,2$. The question is now what will happen for model A
when $d\ge3$ to the (gradient) Gibbs measure, that is, known
to exist without disorder, once we allow for a random environment?

For model B, one can reason similarly as for $d=1,2$ in model A (see
Theorem~1.1 in~\cite{KO}) to show
that there exists no infinite-volume random Gibbs measure if $d=1,2$.
We are interested here in the question whether there exists a random
infinite-volume gradient Gibbs measure for $d\ge1,2$.

To give an intuitive idea of what we can expect, we look next in some
detail at model A in the
special case of a Gaussian (gradient) Gibbs measure where $V (s)
=s^2/2$. In this case one can do explicit computations, and for any
fixed configuration $\xi$, the
finite-volume Gibbs measure with zero boundary condition $\nu_\L
^0[\xi
]$ has expected value
\[
\int\nu_\L^0[\xi](\ormd\phi)(\phi(x))=\sum_{z\in\L}G_\L
(x,z)\xi
(z) \qquad  \mbox{for every fixed } x\in\L,
\]
where $G_\L(x,y)$ denotes the Green's function (see Section~\ref{greenfunction} below for
a~rigorous definition). Due to the properties of the Green's function, the
right-hand side of the equation above diverges as $|\L|\rightarrow
\infty
$ for $d=3,4$ by the Kolmogorov\vadjust{\goodbreak} three series theorem. This hints to the
nonexistence in $d=3,4$ of the infinite-volume $\phi$-Gibbs measure,
which is proved in the \hyperref[app]{Appendix} for the Gaussian case. For the
corresponding gradient Gibbs measure $\mu_\L^0[\xi]$, the expected value
\begin{eqnarray}
\int\mu_\L^0[\xi](\ormd\eta)\bigl(\phi(x)-\phi(y)\bigr)=\sum_{z\in\L
}\bigl(G_\L
(x,z)-G_\L(y,z)\bigr)\xi(z)
\nonumber
\\[-8pt]
\eqntext{\mbox{for every fixed } (x,y)\in\bzd\cap
(\Lambda
\times\Lambda),}
\end{eqnarray}
converges as $|\L|\rightarrow\infty$ for $d\ge3$ and diverges for
$d=1,2$. Coupled with standard tightness arguments, this convergence
for $d\ge3$ gives the existence of the infinite-volume gradient Gibbs
measure in the Gaussian case.
%
%

The main result of our paper, on the existence of shift-covariant
gradient Gibbs measures with given tilt $u\in\R^d$, is the following:
\begin{thmm}
\label{thm1}
\textup{(a)} \textup{(Model A)} Let $d\ge3$, $\E(\xi(0))=0$ and
$u\in\R
^d$. Assume that $V$ satisfies (\ref{tag2}) and (\ref{tag22}). Then
there exists at least one shift-covariant random gradient Gibbs measure
$\xi\rightarrow\mu[\xi]$ with tilt $u$, that is, with
%
%
\begin{eqnarray}
\label{tilt1}
\E\biggl(\int\mu[\xi](\ormd\eta)\eta(b)\biggr)=\langle u,y_b-x_b
\rangle
\nonumber
\\[-8pt]
\\[-8pt]
\eqntext{\mbox{for all bonds } b=(x_b,y_b)\in(\Z^d)^*.}
\end{eqnarray}
Moreover $\mu[\xi]$ satisfies the integrability condition
%
%
\begin{equation}
\label{intcond}
\E\int\mu[\xi](\ormd\eta)(\eta(b))^2<\infty \qquad\mbox{ for all
bonds } b\in
(\zd)^*.
\end{equation}
\begin{enumerate}[(a)]
\item[(b)] \textup{(Model B)} Let $d\ge1$ and $u\in\R^d$. Assume that
$V$ satisfies (\ref{tag23}). Then there exists at least one
shift-covariant random gradient Gibbs measure $\omega\rightarrow\mu
[\omega]$ with tilt $u$, that is, with
%
%
\begin{eqnarray}
\label{tilt2}
\E\biggl(\int\mu[\omega](\ormd\eta)\eta(b)\biggr)=\langle u,y_b-x_b
\rangle
\nonumber
\\[-8pt]
\\[-8pt]
\eqntext{\mbox{for all bonds } b=(x_b,y_b)\in(\Z^d)^*.}
\end{eqnarray}
Moreover $\mu[\omega]$ satisfies the integrability condition
%
%
\begin{equation}
\E\int\mu[\omega](\ormd\eta)(\eta(b))^2<\infty\qquad \mbox{for all
bonds } b\in(\zd)^*.
\end{equation}
\end{enumerate}
\end{thmm}

For model A we also show by similar arguments as in~\cite{EK} the following:
\begin{thmm}[(Model A)]
\label{non-exist of Gibbs}
Let $d\ge3$ and assume that $\E(\xi(0))\neq0$.
Then there exists no shift-covariant gradient Gibbs measure $\mu[\xi
]$ with
\[
\E\biggl|\int\mu[\xi](\mathrm{d}\eta)V'(\eta(b))\biggr|<\infty \qquad\mbox{for all
bonds } b=(x,y)\in(\zd)^*.
\]
\end{thmm}

The techniques used to prove existence in the nonrandom continuous
interface model are based on the Brascamp--Lieb inequality and on
shift-invariance, which techniques do not work in our random settings;
the lack of shift-invariance in our models means that the
Brascamp--Lieb inequality is not enough to ensure tightness of the
finite-volume gradient Gibbs measures $(\mu_\L^\rho[\xi])$,
respectively, of $(\mu_\L^\rho[\omega])$, as is the case in
the model without disorder (see the \hyperref[app]{Appendix} for a more detailed
explanation of the Brascamp--Lieb inequality and why it fails in the
case of our models in a disordered setting).
We will prove the existence result for model A and sketch it for model
B. To prove our result for model A, we are using surface tension bounds
to establish tightness of a sequence of spatially averaged
finite-volume gradient Gibbs measures for each realization of the
disorder, whose limit along a~deterministic subsequence we extract
(using a result in~\cite{KOM}) and we prove that it is a
shift-covariant random gradient Gibbs measure.

To complement our analysis of the two models, we will also investigate
under what assumptions on the disorder $\xi$, respectively, on
$V^\omega
_{(x,y)}$, the surface tension $\sigma[\xi](u)$, respectively, $\sigma
[\omega](u)$, exists and under what assumptions it does not exist.
Moreover we will prove that when it exists, the surface tension is $\P
$- a.s. independent of the disorder. The surface tension bounds
established in Theorem~\ref{non-exist}(b) are used later to prove
tightness of the finite-volume spatially averaged Gibbs measures,
averaged over the disorder. To state our surface tension result, let
$a,l\in\Z^d,a=(a_1,\ldots,a_d), l=(l_1,\ldots,l_d)$, with $a_i<l_i,
i=1,2,\ldots,d$, and let
%
%
\begin{equation}
\label{lal}
\bar\Lambda^{a,l}:=\{z\in\Z^d\dvtx a_i\le z_i\le l_i \mbox{ for
all } i=1,2,\ldots,d\}.
\end{equation}
For any $n\in\Z$, we denote by $a+n:=(a_1+n,\ldots,a_d+n)$ and
by $an:=(a_1n,\ldots,a_dn)$. In view of Theorem~\ref{non-exist}(a) and of
Remark~\ref{surf12} below, we have

\begin{thmm}[(Model A)]\label{nonexist1}
The infinite-volume surface tension does not exist
if $d=1,2$ or if $d\ge3$ and $\E(\xi(0))\neq0$.
\end{thmm}

For $d\ge3$ and $\E(\xi(0))= 0$, we prove
\begin{thmm}
\label{subadit}
\textup{(1)} \textup{(Model A)}
Let $d\ge3$ and assume that $\E(\xi(0))=0$ and $u\in\R^d$. Then if $V$
satisfies (\ref{tag2}) and (\ref{tag22}), we have:
\begin{longlist}[(a)]
\item[(a)] $\sigma[\xi](u):=\lim_{ N\rightarrow\infty}\sigma
_{\Lambda
_{N}}[\xi](u)$ exists for $\P$-almost all $\xi$ and in $L^1$ and
\[
{\sigma}[\xi](u)=\lim_{n\rightarrow\infty}\frac{1}{n^d}\lim
_{m\rightarrow\infty}\frac{1}{m^d}\sum_{a_i\in\N, 1\le a_i\le
m,i=1,\ldots, d}\sigma_{\bar\Lambda^{(a-1)n,an}}(u),
\]
where the limits in $m\rightarrow\infty$ and in $n\rightarrow\infty$
are in $L^1$.
%
\item[(b)] $\sigma[\xi](u)$ is independent of $\xi$, with
\[
\sigma[\xi](u)=\lim_{ N\rightarrow\infty}\E(\sigma_{\Lambda
_{N}}[\xi](u))=:\sigma(u)\qquad \mbox{for } \P\mbox{-almost all } \xi.\vadjust{\goodbreak}
\]
\end{longlist}
\begin{enumerate}[(2)]
\item[(2)] \textup{(Model B)} Let $d\ge1$. Then $\sigma[\omega](u)$
satisfies \textup{(a)--(b)} above, with $\omega$ replacing~$\xi$ in the results.
\end{enumerate}
\end{thmm}

The presence of the disorder and of the Green's functions make the question
of existence of the surface tension more delicate to handle than in the
nonrandom case, where the answer is fairly straightforward. In order to
prove existence of the surface tension for our disordered system, we
prove (almost)-subadditivity of the finite-volume surface tension, in
order to apply ergodic theorems for subadditive processes.

A natural question to ask is whether in our disordered models a random
gradient Gibbs measure is uniquely determined by the tilt as in the
nonrandom settings of~\cite{CD} or~\cite{FS}. This is work in progress
by the same authors and will be addressed in a future paper.

The rest of the paper is organized as follows: In Section~\ref{sec2} we recall
the definition and some basic properties of the Green's function and we
prove a~strong law of large numbers (SLLN) involving the Green's function,
which are necessary for the proof of our main Theorem~\ref{thm1} and
for the surface tension results; we also recall in Section~\ref{sec2} two
subadditivity propositions used for the proof of the surface tension
existence. In Section~\ref{sec3}, we study model A. In Section~\ref{tsf}, we prove
Theorem~\ref{non-exist}, and, respectively, Theorem~\ref{subadit}, for
nonexistence and, respectively, for existence of the surface tension.
In Section~\ref{esvg1}, we formulate and prove Theorem~\ref{thm1}, our main
result on the existence of shift-covariant random gradient Gibbs
measures. Section~\ref{sec4} deals with the corresponding results for model B.
Finally, the \hyperref[app]{Appendix} explains why the infinite-volume Gibbs measure
for model A does not exist for \mbox{$d=3,4,$} and provides a more detailed
explanation of the Brascamp--Lieb inequality.

\section{Preliminary notions}\label{sec2}

\subsection{Green functions on $\Z^d$}
\label{greenfunction}
We first review a few facts about Green's functions.

Let $A$ be an arbitrary subset in $\Z^d$ and let $x\in A$ be fixed. Let
$\P_x$ and~$\E_x$ be the probability law and expectation, respectively,
of a simple random walk $X:=(X_k)_{k\geq0}$ starting from $x\in\Z^d$;
Green's function $G_A(x,y)$ is the expected number of visits to
$y\in A$ of the walk $X$ killed as it exits $A$, that is,
\[
G_A(x,y)=\E_x\Biggl[\sum_{k=0}^{\tau_A-1} 1_{(X_k=y)}\Biggr]=\sum
_{k=0}^\infty\P_x(X_k=y,k<\tau_A),\qquad  y\in\Z^d,
\]
where $\tau_A=\inf\{k\ge0\dvtx X_k\in A^c\}$. We will state first some
well-known properties of the Green's functions. To avoid exceptional cases
when $x=0$, let us denote by $]|x|[\,=\max\{|x|,1\}$, where $|x|$ is the
Euclidian norm.
\begin{prop}
\label{properties}
\begin{enumerate}[(iii)]
\item[(i)] If $d\ge3$, then $\lim_{N\rightarrow\infty}G_{\Lambda
_N}(x,y):=G(x,y)$ exists for all $x,y\in\Z^d$ and as
$|x-y|\rightarrow
\infty$,
\[
G(x,y)=\frac{a_d}{|x-y|^{d-2}}+O(|x-y|^{1-d}),
\]
with $a_d=\frac{2}{(d-2)w_d}$, where $w_d$ is the volume of the unit
ball in $\R^d$.
\item[(ii)] Let $B_r=\{x\in\Z^d\dvtx |x|<r\}$; then for $x\in B_N$
\[
G_{B_N}(0,x)=
\cases{
\displaystyle\frac{2}{\pi}\log\frac{N}{]|x|[}+o\biggl(\frac{1}{]|x|[}\biggr)+O
\biggl(\frac{1}{N}\biggr),&\quad $\mbox{if } d=2,$\vspace*{2pt}\cr
\displaystyle\frac{2}{(d-2)w_d}\bigl[ ]|x|[^{2-d}-N^{2-d}+O( ]|x|[^{1-d}
)\bigr], & \quad $\mbox{if } d\ge3.$}
\]
Let $\varepsilon>0$. If $x\in B_{(1-\varepsilon)N}$, the following
inequalities hold:
\[
G_{B_{{\varepsilon N}}}(0,0)\le G_{B_N}(x,x)\le G_{B_{2N}}(0,0).
\]
\item[(iii)] $G_A(x,y)=G_A(y,x)$.
\item[(iv)] $G_A(x,y)\le G_B(x,y)$, if $A\subset B$.
\item[(v)] If $x\in B_N$, then
\[
N^2-|x|^2\le\E_x(\tau_{B_N})\le(N+1)^2-|x|^2.
\]
\end{enumerate}
\end{prop}

For proofs of (i), (iii) and (iv) from Proposition \ref
{properties} above we refer to Chapter 1 from~\cite{Lawl1}, for proof
of (ii) we refer to Lemma 1 from~\cite{L} and for proof of (v) we refer
to Lemma 2 from~\cite{L}.

The result we state next will be used to prove Theorem~\ref{non-exist}.
\begin{prop}
\label{green}
There exists $N_0$ sufficiently large such that for all $N\ge N_0$, we have
\[
\frac{d+1}{d+2}w_dN^{2}(N-1)^d\le\sum_{x,y\in\L_N}G_{\Lambda
_N}(x,y)\le
\bigl(N\sqrt{d}\bigr)^d dw_d\biggl[(N+1)^2-\frac{N^2}{d+2}\biggr].
\]
\end{prop}

\begin{pf}
Note first that since $G_{B_N}$ is symmetric, we have
%
%
\begin{equation}
\label{taueqn}
\E_x(\tau_{B_N})=\sum_{y\in B_{N}} G_{B_N}(x,y)=\sum_{y\in B_{N}}
G_{B_N}(y,x).
\end{equation}
\textit{The upper bound}:
Using Proposition~\ref{properties}(iv) for the first inequality, (\ref
{taueqn}) for the second inequality and Proposition~\ref{properties}(v)
for the third inequality, we have for $N$ large enough
\begin{eqnarray*}
\sum_{x,y\in\L_N}G_{\Lambda_N}(x,y)&\le&\sum_{x,y\in B_{N\sqrt
{d}}}G_{B_{N\sqrt{d}}}(x,y)=\sum_{x\in B_{N\sqrt{d}}}\E_x(\tau
_{B_{N\sqrt{d}}})\\
&\le&\sum_{x\in B_{N\sqrt{d}}}[(N+1)^2-|x|^2
]\\
&\le&(N+1)^2d \bigl(N\sqrt{d}\bigr)^d w_d -w_d\int_0^{N\sqrt
{d}}r^{d+1}\,\ormd r\\
&=&\bigl(N\sqrt{d}\bigr)^d dw_d\biggl[(N+1)^2-\frac
{N^2}{d+2}\biggr].
\end{eqnarray*}
\textit{The lower bound}: We have $B_{N}\subset\L_N$. Then by using
Proposition~\ref{properties}(iv), (v) and~(\ref{taueqn}), we have for
$N$ large enough
\begin{eqnarray*}
\sum_{x,y\in\L_N}G_{\Lambda_N}(x,y)&\ge&\sum_{x,y\in
B_{N}}G_{B_{N}}(x,y)\ge\sum_{x\in B_{N}}[N^2-|x|^2]\\
&\ge&
N^{2}(N-1)^d w_d -w_d\int_0^{N-1}r^{d+1}\,\ormd r\\
&\ge&\frac{d+1}{d+2}w_d(N-1)^{d}N^2.
\end{eqnarray*}
\upqed\end{pf}

We will use the next result in the proof of Proposition~\ref{norafone}.
\begin{prop} \label{monoform}
Let $d\geq1$ and let $\L_1\subset\L_2 \subset\Z^d$. Then we have for
all $\x\in\R^{\L_2}$
%
%
\begin{equation}
\langle\xi, G_{\L_1} \xi\rangle_{\L_1}\leq\langle\xi, G_{\L
_2} \xi
\rangle_{\L_2},
\end{equation}
where $\langle\xi,G_\L\xi\rangle_\L:=\sum_{x,y\in\L}\xi
(x)G_\L(x,y)\xi
(y)$ and where $G_{\L}:= (G_{\L}(x, y))_{x,y\in\L}$.
\end{prop}

\begin{pf} A proof of this statement can be found, for
example, in~\cite{Sim}.~%
\end{pf}

\subsection{Strong law of large numbers}
We will need the following strong law of large numbers (SLLN) in the
proof of Theorems~\ref{non-exist} and~\ref{subadit}.
\begin{prop}
\label{sllngree}
Let $(\xi(x))_{x\in\Z^d}$ be i.i.d. with $\E(\xi^2(0))<\infty$.
For all
$d\ge3$, we have
%
%
\begin{equation}
\label{sllngreen1}\qquad
\lim_{N\rightarrow\infty}\frac{\langle\xi,G_{\L_N}\xi\rangle
_{\Lambda
_N} -\sum_{x,y\in\Lambda_N}\E(\xi(x)\xi(y)) G_{\L
_N}(x,y)}{N^{d}}=0 \qquad \mbox{a.s.}
\end{equation}
\end{prop}

\begin{pf} Let the variance w.r.t. $\P$ be denoted by
$\bvar
$ and let
%
\begin{eqnarray*}
S_N&:=&\frac{\sum_{x,y\in\Lambda_N\atop x\neq y}[\xi(x)-\E(\xi
(x))][ \xi(y)-\E(\xi(y))]G_{\L_N}(x,y)}{N^{d}},
\\
%
%
S_N'&:=&\frac{\sum_{x,y\in\Lambda_N\atop x\neq y}[\xi(x)-\E(\xi
(x))]\E(\xi(y))G_{\L_N}(x,y)}{N^{d}}\quad \mbox{and}\\
R_N&:=&\frac{\sum
_{x\in\Lambda_N}[\xi^2(x)-\E(\xi^2(x))]G_{\L_N}(x,x)}{N^{d}}.
\end{eqnarray*}
%
Note that proving (\ref{sllngreen1}) is the same as proving that
\[
\lim_{N\rightarrow\infty} S_N=0 ,\qquad  \lim_{N\rightarrow\infty}
S_N'=0\quad \mbox{and}\quad \lim_{N\rightarrow\infty} R_N=0 \qquad\mbox{a.s.}
\]
Using the independence of the $(\xi(x))_{x\in\zd}$ for the equality
below, Proposition~\ref{properties}(iv) for the first inequality below
and (ii) for the second one, we have
%
%
\begin{eqnarray}
\label{bc3}
\E(S_N^2)&=&\frac{\bvar^2(\xi)}{N^{2d}}\sum_{x,y\in\Lambda
_N,x\neq
y}G^2_{\L_N}(x,y)\nonumber\\
&\le&\frac{\bvar^2(\xi)}{N^{2d}}\sum_{x,y\in
B_{N\sqrt
{d}},y\neq x} G^2_{B_{N\sqrt{d}}}(x,y)\nonumber\\
&\le&\frac{\bvar^2(\xi)}{N^{2d}}\biggl(\frac{2}{(d-2)w_d}\biggr)^2 \sum
_{x,y\in B_{N\sqrt{d}},y\neq x}\biggl(\frac{1}{|x-y|^{2d-4}}+O(1)
\biggr)\\
&\le& C(w_d)\frac{\bvar^2(\xi)}{N^{2d}}\sum_{x\in B_{N\sqrt{d}}}
\biggl(\int_1^{N\sqrt{d}}\frac{1}{r^{d-3}}\,\mathrm{d}r+O(1)\biggr)\nonumber\\
&\le&\tilde
{C}(w_d,d)\frac{\bvar^2(\xi)}{N^{d-1}}.\nonumber
\end{eqnarray}
Fix $\varepsilon>0$. By means of (\ref{bc3}), we get
\[
\sum_{N=1}^{\infty}\P(|S_N|\ge\varepsilon)\le\tilde
{C}(w_d,d)\frac{{\bvar}^2(\xi^2)}{\varepsilon^2}\sum_{N=1}^\infty
\frac{1}{
N^{d-1}}<\infty
\]
and therefore by Borel--Cantelli
\[
\lim_{N\rightarrow\infty}\sup|S_N|\le\varepsilon\qquad \mbox{a.s.}, \quad\mbox{from which}\quad
\lim_{N\rightarrow\infty}S_N=0\qquad \mbox{a.s.}
\]
The proof that $\lim_{N\rightarrow\infty} S'_N=0 \mbox{ a.s.}$ follows
the same pattern as the proof for~$S_N$, and will be omitted. We will
proceed next with the proof of\break
$\lim_{N\rightarrow\infty} R_N=0 \mbox{ a.s.}$ Let $\varepsilon>0$ be
arbitrarily fixed and denote for simplicity of notation $\tau(x):=(\xi
^2(x)-\E(\xi^2(x)))$. Take $M=M(\varepsilon)>0$ such that $\E(|\tau
(x)|1_{|\tau(x)|>M})\le\varepsilon$ and define
\[
R'_N=\frac{\sum_{x\in\Lambda_N}G_{\Lambda_N}(x,x)[\tau(x)1_{|\tau
(x)|>M}-\E(\tau(x)1_{|\tau(x)|>M})]}{N^d}
\]
and
\[
R''_N=\frac{\sum_{x\in\Lambda_N}G_{\Lambda_N}(x,x)[\tau
(x)1_{|\tau
(x)|\le M}-\E(\tau(x)1_{|\tau(x)|\le M})]}{N^d}.
\]
Using Proposition~\ref{properties}(ii) and (iv) to find $C>0$ such that
$|G_{\Lambda_N}(x,x)|\le C$, uniformly in $N$ and $x\in\Lambda_N$, and
using the SLLN for i.i.d. random variables with finite first moment, we get
\begin{eqnarray*}
|R_N'|&\le& C\frac{\sum_{x\in\Lambda_N}[|\tau(x)|1_{|\tau
(x)|>M}+\E
(|\tau(x)|1_{|\tau(x)|>M})]}{N^d}\\
&\le&2C\E\bigl(|\tau
|1_{|\tau|>M}\bigr)\bigl(1+o(1)\bigr)\\
&\le& 2C\varepsilon\bigl(1+o(1)\bigr).
\end{eqnarray*}
Therefore
%
%
\begin{equation}
\label{sup1}
\limsup_{N\rightarrow\infty}|R_N'|\le2C\varepsilon\qquad \mbox{a.s.},
\end{equation}
from which we get $R_N'\rightarrow0$ a.s. Since the summands in
$R_N''$ are uniformly bounded and independent, by a standard fourth
moment bound, Markov inequality and Borel--Cantelli, we have $R_N''
\rightarrow0$ a.s. This concludes the proof of the proposition.
\end{pf}

\subsection{Ergodic theorems for multiparameter subadditive processes}
$\!\!\!$For \mbox{$N\in\N$}, let $\Lambda_{[0,N]}:=[0,N]^d\cap\Z^d$, let $\Z
^d_+:=\{
z\in\Z^d\dvtx 0\le z_i \mbox{ for all } i=1,2,\ldots,d\}$ and let
\begin{eqnarray*}
\mathcal{A}&:=&\{\Lambda\subset\Z^d_+\dvtx \Lambda=\bar\Lambda
^{a,l} \mbox{ for
some } a,l\in\Z^d_+, a=(a_i)_{1\le i\le d},l=(l_i)_{1\le i\le
d}, \\
&&\hspace*{202pt}\phantom{\{}\mbox{with }  a_i<l_i,1\le i\le d\},
\end{eqnarray*}
where we recall that $\bar\Lambda^{a,l}$ was defined in (\ref{lal}).
For any finite set $\Lambda\in\Z^d$ and for any $z\in\Z^d$, we denote
$\Lambda+z:=\{x+z\dvtx x\in\Lambda\}$.

We will use the two propositions below to prove a.s. and $L^1$
convergence of the surface tension. The first proposition is an ergodic
theorem for superadditive processes from~\cite{AKKRE}:
\begin{prop}
\label{akk}
Let $(\tau_z)_{z\in\Z^d_+}$ be a measurable semigroup of
measure-preserving transformations on $(\Omega,\mathcal{F}, \P)$. Let
$(W_I)_{I\in\mathcal{A}}$ be a family of real-valued random variables on
$(\Omega,\mathcal{F}, \P)$ such that a.s.:
\begin{longlist}[(a)]
\item[(a)] $W_I\circ\tau_z=W_{I+z}$.
\item[(b)] (The subadditivity condition) If $\bigcup_{i=1}^n I_i=I\in
\mathcal{A}$ with $(I_i)_{i=1,2,\ldots,n}$ pairwise disjoint in $\mathcal{A}$,
then $W_I\le\sum_{i=1}^n W_{I_i}$.
\item[(c)]
\begin{eqnarray}
\inf|I|^{-1}\int W_I \,\ormd\P>-\infty
\nonumber\\[-3pt]
\eqntext{\mbox{the infimum being taken
over all }  I\in\mathcal{A} \mbox{ with }  |I|>0,}
\end{eqnarray}
where $|I|$ denotes the cardinality of the finite set $I$.
\end{longlist}
Then $\lim_{N\rightarrow\infty}N^{-d} W_{\Lambda_{[0,N]}} \mbox
{ exists a.s.}$\vadjust{\goodbreak}
\end{prop}

The second proposition is Theorem 2.1 from~\cite{Schur}. In what
follows, $x^+$ denotes the positive part of $x\in\R$.
\begin{prop}
\label{schurl1}
Let $(W_I)_{I\in\mathcal{A}}$ be a family of real-valued random variables
on $(\Omega,\mathcal{F}, \P)$ such that:
\begin{longlist}[(a)]
\item[(a)] If $\bigcup_{i=1}^n I_i=I\in\mathcal{A}$ with
$(I_i)_{i=1,2,\ldots,n}$ pairwise disjoint in $\mathcal{A}$, then $\E
(W_I-\sum_{i=1}^n W_{I_i})\le0.$
\item[(b)] $E(W_{I+z})=E(W_I)$ for all $I\in\mathcal{A}$ and $z\in\Z^d_+$.
\item[(c)] $E(W_{I+z}^+)=E(W_I^+)$ for all $I\in\mathcal{A}$ and $z\in
\Z^d_+$.
\item[(d)]
\begin{eqnarray}
\inf|I|^{-1}\int W_I \,\ormd\P>-\infty
\nonumber\\[-3pt]
\eqntext{\mbox{the infimum being taken
over all }  I\in\mathcal{A} \mbox{ with }  |I|>0.}
\end{eqnarray}
\item[(e)] Assume that for every $a,l\in\Z^d_+, a=(a_i)_{1\le i\le
d},l=(l_i)_{1\le i\le d}$, the collection of random variables $
(W_{a,l})_{a,l\in\Z^d_+}$, with $W_{a,l}:=W_{\bar\Lambda
^{(a-1)n,an}},$ is stationary with respect to all translations in $\Z
^d$ of form $(a,l)\rightarrow(a+v,l+v)$.
\end{longlist}
Then
\[
\lim_{N\rightarrow\infty}N^{-d} W_{\Lambda_{[0,N]}}=W_{\infty
}\qquad \mbox{exists in } L^1,
\]
where
\[
W_{\infty}=\lim_{n\rightarrow\infty}\frac{1}{n^d}\lim
_{m\rightarrow
\infty}\frac{1}{m^d}\sum_{1\le a_i\le m,i=1,\ldots, d}W_{\bar
\Lambda
^{(a-1)n,an}}
\]
and where the limits in $m\rightarrow\infty$ and in $n\rightarrow
\infty
$ are in $L^1$.
\end{prop}

Both Proposition~\ref{akk} and Proposition~\ref{schurl1} can be stated
and proved for sets $\bar{\mathcal{A}}$ in~$\Z^d$ of form
\begin{eqnarray*}
\bar{\mathcal{A}}&:=&\{\Lambda\subset\Z^d\dvtx \Lambda=\bar\Lambda
^{a,l} \mbox{ for some } a,l\in\Z^d, a=(a_i)_{1\le i\le d},l=(l_i)_{1\le i\le
d},\\
 &&\hspace*{199pt}\phantom{\{}\mbox{with } a_i<l_i,1\le i\le d\},
\end{eqnarray*}
instead of just for sets $\mathcal{A}$ in $\Z^d_+$.

\section{Model A}\label{sec3}
This section is structured as follows: in Section~\ref{tsf1} we
prove Theorem~\ref{non-exist}, on the nonexistence of the surface
tension when $\E(\xi(0))\neq0$; in Section~\ref{tsf2} we prove
Theorem~\ref{subadit}, on the existence of the surface tension when
$d\ge3$ and $\E(\xi(0))=0$, by means of subadditivity arguments. In
Section~\ref{esvg1} we prove Proposition~\ref{averaged}, on the
tightness of the finite-volume gradient Gibbs measures $(\mu_\L^\rho
[\xi])_{\L\in\Z^d}$ averaged over the disorder, from which we
derive the existence of the random infinite-volume gradient Gibbs
measure averaged over the disorder. This tightness result is
instrumental in Section~\ref{esvg2}, in our proof of existence of
the infinite-volume random gradient Gibbs measure.

\subsection{The surface tension}
\label{tsf}


\subsubsection{\texorpdfstring{Nonexistence of the surface tension when $\E(\xi
(0))\neq0$}{Nonexistence of the surface tension when E(xi(0)) /= 0}}
\label{tsf1}

We prove in this subsection that the surface tension does not exist
when $\E(\xi(0))\neq0$, and when $\E(\xi(0))=0$ we give upper and
lower bounds on $\sigma_{\Lambda_N}[\xi](u)$, uniformly in $\Lambda_N$.
\begin{thmm}
\label{non-exist}

Let $d\ge3$. Assume that $V$ satisfies (\ref{tag2}) and (\ref{tag22}).
Recall that $(\xi(x))_{x\in\Z^d}$ are i.i.d. with finite second moments.
\begin{enumerate}[(a)]
\item[(a)] If $\E(\xi(0))\neq0$, then for all $u\in\R^d$
\[
S_1\le\liminf_{N\rightarrow\infty}\frac{\sigma_{\Lambda_N}[\xi
](u)}{N^{2}}\le\limsup_{N\rightarrow\infty}\frac{\sigma_{\Lambda
_N}[\xi
](u)}{N^{2}}\le S_2\qquad \mbox{for } \P\mbox{-almost all } \xi,
\]
where
\[
S_1:=-\frac{w_d}{2A(d+2)}\E^2(\xi(0))\quad \mbox{and}\quad S_2:=-\frac
{w_d(d+1)}{4C_2(d+2)} \bigl(\sqrt{d}\bigr)^d\E^2(\xi(0)).
\]
\item[(b)] If $\E(\xi(0))=0$, then
%
%
\begin{eqnarray}
\bar{S}_1(u) \le\liminf_{N\rightarrow\infty}\sigma_{\Lambda
_N}[\xi
](u)\le\limsup_{N\rightarrow\infty}\sigma_{\Lambda_N}[\xi](u)\le
\bar
{S}_2(u)
\nonumber
\\[-8pt]
\\[-8pt]
\eqntext{\mbox{for } \P\mbox{-almost all } \xi,}
\end{eqnarray}
where
\begin{eqnarray*}
\bar{S}_1(u)&:=&\sigma^{A}[\xi=0](u=0)-\frac{w_d}{A(d-2)}\E(\xi
^2(0))+A(1+|u|^2)-2dB,
\\
\bar{S}_2(u)&:=& \sigma^{C_2/2}[\xi=0](u=0)-\frac{w_d}{2C_2(d-2)}\E
(\xi
^2(0))+\frac{C_2}{2}(1+|u|^2)\\
&&{}+2dV(0).
\end{eqnarray*}
For a $C>0$, we defined by $\sigma_{\L}^C[\xi=0](u=0)$ and $\sigma
^{C}[\xi=0](u=0)$ the finite-volume and infinite-volume surface
tensions corresponding to model A without disorder, with potential
$V(x)=Cx^2$ and tilt $u=0$.
\end{enumerate}
\end{thmm}

In particular, the above theorem shows that if $\E(\xi
(0))\neq0$, then the surface tension does not exist as the
finite-volume surface tension $\log Z^{\psi_u}_{\L_N}[\xi]$ is of order
$N^{d+2},$ and not of order $N^{d},$ as would normally be expected (and
as indeed is the case in the nondisordered case). The reason that the
$N^{d+2}$ exponent comes up is mainly due to the appearance of
the Green's function in the formulas for the upper/lower bounds for the
finite-volume surface tension. When $\E(\xi(0))\neq0$, the~terms in
the upper/lower bounds involve
double sums over the Green's function
of the form $\sum_{x,y \in\Lambda_N} G_{\Lambda_N}(x,y)$, which are of
order $N^{d+2}$.\vadjust{\goodbreak}

\begin{pf*}{Proof of Theorem \protect\ref{non-exist}}
We will use the bounds for $V$ from (\ref{tag2}) and (\ref{tag22}) to
obtain upper and lower bounds for $\sigma_{\Lambda_N}[\xi]$ in terms of
surface tensions for the nondisordered model with
quadratic potentials. The claims in (a) and (b) will follow then easily
by an application of Proposition~\ref{sllngree}. The explicit
computations follow below.

We will start by proving a lower bound for $\sigma_{\Lambda_N}[\xi
](u)$. As $V(s)\ge As^2-B$, we get from (\ref{surftensdef})
%
%
\begin{eqnarray}
\label{lowerbound}
&&\sigma_{\Lambda_N}[\xi](u)\nonumber\\
&&\qquad\ge -\frac{1}{2|\Lambda_N|} \mathop{\sum_{x,y\in\L_N\cup\partial\L
_N}}_{|x-y|=1} B\nonumber\\
&&\qquad\hspace*{9pt}{}-\frac{1}{|\Lambda_N|}\log\int\exp\biggl(
-\frac{A}{2} \mathop{\sum_{x,y\in\L_N}}_{|x-y|=1}\bigl(\phi(x)-\phi(y)\bigr)^2
\nonumber\\
&&\hspace*{118pt}{}-A \mathop{\sum_{x\in\L_N,y\in\partial\L_N}}_{|x-y|=1}\bigl(\phi
(x)-{{{\psi}}_u}(y)\bigr)^2\nonumber\\
&&\hspace*{189pt}{} +\sum_{x\in\L_N}\xi(x)
\phi(x)
\biggr) \,\ormd\phi_{\L_N}\\
&&\qquad=-2dB-\frac{\sum_{x\in\L_N}\xi(x) (x\cdot u)}{|\Lambda
_N|}\nonumber\\
&&\qquad\hspace*{9pt}{}-\frac
{1}{|\L_N|}\log\int\exp\biggl(
-\frac{A}{2}\mathop{\sum_{x,y\in\L_N}}_{|x-y|=1}\bigl(\tilde\phi(x)-\tilde
\phi
(y)+ (x-y)\cdot u\bigr)^2 \nonumber\\
&&\hspace*{118pt}{}-A \mathop{\sum_{x\in\L
_N,y\in\partial\L_N}}_{|x-y|=1}\bigl(\tilde\phi(x)+(x-y)\cdot u\bigr)^2
\nonumber\\
&&\hspace*{220pt}{}+\sum_{x\in\L_N}\xi(x) \tilde\phi(x)
\biggr)\,\ormd\tilde\phi_{\L_N},\nonumber
\end{eqnarray}
where for the equality we used the change of variables $\phi(x)=\tilde
\phi(x)+x\cdot u$ for all $x\in\L_N$. To simplify (\ref
{lowerbound}) we
will show next that
\begin{eqnarray}
\label{minquad}
\qquad&&\frac{1}{2}\mathop{\sum_{x,y\in\L_N}}_{|x-y|=1}\bigl(\tilde\phi(x)-\tilde
\phi
(y)+(x-y)\cdot u\bigr)^2+\mathop{\sum_{x\in\L_N,y\in\partial\L_N}}_{|x-y|=1}\bigl(\tilde\phi(x)+(x-y)\cdot u\bigr)^2\nonumber\\[-1pt]
&&\qquad=\frac{1}{2}\mathop{\sum_{x,y\in\L_N}}_{|x-y|=1}\bigl[\bigl(\tilde\phi
(x)-\tilde\phi(y)\bigr)^2+\bigl((x-y)\cdot u\bigr)^2\bigr]\\[-1pt]
&&\qquad\quad{}+\mathop{\sum_{x\in\L_N,y\in
\partial\L_N}}_{|x-y|=1}\bigl[(\tilde\phi(x))^2+\bigl((x-y)\cdot
u\bigr)^2\bigr].\nonumber
\end{eqnarray}
By expanding the square, (\ref{minquad}) follows from
\[
\mathop{\sum_{x,y\in\L_N}}_{|x-y|=1}[\tilde\phi(x)-\tilde\phi
(y)][(x-y)\cdot u]+2\mathop{\sum_{x\in\L_N,y\in\partial\L
_N}}_{|x-y|=1}\tilde\phi(x)[(x-y)\cdot u]=0,
\]
which can be easily seen to be true by summing over
bonds along lines in each coordinate direction. Plugging the identity
from (\ref{minquad}) into (\ref{lowerbound}), we get
%
%
\begin{eqnarray}
\label{1a}
&&\sigma_{\Lambda_N}[\xi](u)\nonumber\\[-1pt]
%
&&\qquad\ge-2dB+\frac{A}{2|\L_N|}\mathop{\sum_{{x,y\in\L_N}}}_{|x-y|=1}\bigl((x-y)\cdot u\bigr)^2
\nonumber\\[-1pt]
&&\qquad\quad{}+\frac{A}{|\L_N|}\sum_{{x\in\L,y\in\partial\L_N}\atop
|x-y|=1}\bigl((x-y)\cdot u\bigr)^2-\frac{\sum_{x\in\L_N}\xi(x) (x\cdot
u)}{|\L
_N|}\\[-1pt]
&&\qquad\quad{}-\frac{1}{|\L_N|}\log\int\exp\biggl( -\frac{A}{2}\mathop{\sum_{x,y\in\L
_N}}_{|x-y|=1}\bigl(\tilde\phi(x)-\tilde\phi(y)\bigr)^2\nonumber\\[-1pt]
&&\qquad\hspace*{97pt}{} -A \mathop{\sum_{x\in
\L_N,y\in\partial\L_N}}_{|x-y|=1}(\tilde\phi(x))^2+\sum_{x\in
\L
_N}\xi(x) \tilde\phi(x)\biggr) \,\ormd\tilde\phi_{\L_N}.\nonumber
\end{eqnarray}
%
To compute the integral in (\ref{1a}) we use standard Gaussian calculus
(see, e.g., Proposition 3.1 part (2) from~\cite{FS}) to show that
%
%
\begin{eqnarray}
\label{mult1}
&&\log\int\exp\biggl( -\frac{A}{2}\mathop{\sum_{x,y\in\L_N}}_{
|x-y|=1}\bigl(\tilde\phi(x)-\tilde\phi(y)\bigr)^2\nonumber \hspace*{-25pt}\\[-1pt]
&&\qquad{}-A\mathop{\sum_{x\in\L_N,y\in
\partial\L_N}}_{|x-y|=1}(\tilde\phi(x))^2+\sum_{x\in\L_N}\xi(x)
\tilde\phi(x)\biggr)\,\ormd\tilde\phi_{\L_N}\nonumber\hspace*{-25pt}\\[-1pt]
&&\qquad=\log\int\exp\biggl( -\frac{A}{2}\mathop{\sum_{x,y\in\L_N}}_{
|x-y|=1}\bigl(\tilde\phi(x)-\tilde\phi(y)\bigr)^2 -A\mathop{\sum_{x\in\L_N,y\in
\partial\L_N}}_{|x-y|=1}(\tilde\phi(x))^2\biggr) \,\ormd\tilde\phi
_{\L_N}\hspace*{-25pt} \\[-1pt]
&&\qquad\quad{}+\frac{\langle\xi, G_{\L_N} \xi\rangle_{\L_N}}{2A}
\nonumber\hspace*{-25pt}\\[-1pt]
&&\qquad=-|\L_N|\s_{\L_N}^A[\xi=0](u=0)+\frac{\langle\xi, G_{\L_N}
\xi\rangle
_{\L_N}}{2A}.\nonumber\hspace*{-25pt}
\end{eqnarray}
Plugging (\ref{mult1}) in (\ref{1a}) gives the lower bound for
$\sigma
_{\Lambda_N}[\xi](u)$.

Due to the assumption $V''\le C_2$, we have by Taylor expansion that
$V(s)\le V(0)+\frac{C_2}{2}s^2$; then by the same reasoning as in the
derivation of the lower bound, we get
%
%
\begin{eqnarray}
\label{upbound}
\sigma_{\Lambda_N}[\xi](u)&\le& 2dV(0) +\s_{\L_N}^ {C_2/2}[\xi=0](u=0)
+\frac{C_2}{4|\L_N|}\mathop{\sum_{x,y\in\L_N}}_{|x-y|=1}\bigl((x-y)\cdot
u\bigr)^2\nonumber\\
&&{}+\frac{C_2}{2|\L_N|}\mathop{\sum_{x\in\L,y\in\partial\L_N}}_{|x-y|=1}\bigl((x-y)\cdot u\bigr)^2-\frac{\sum_{x\in\L_N}\xi(x) (x\cdot
u)}{|\L
_N|}\\
&&{}- \frac{\langle\xi, G_{\L_N}
\xi\rangle_{\L_N}}{4C_2|\L_N|}.\nonumber
\end{eqnarray}
The upper bound follows now from (\ref{upbound}), by noting that for
all $C>0$, $\sigma_{\L}^C[\xi=0](u)\rightarrow\sigma^C[\xi
=0](u)\in
(-\infty,\infty)$ as $|\L|\rightarrow\infty$ (for a proof of this, see
Proposition 1.1 in~\cite{FS}).
\begin{longlist}[(a)]
\item[(a)]
The statement follows now from (\ref{lowerbound}), (\ref{upbound}),
Proposition~\ref{sllngree} and Proposition~\ref{green} by noting that
for very large $N$
\[
w_d\frac{d+1}{d+2}\E^2(\xi(0))\le\frac{1}{N^{d+2}}\E(\langle\x
,(G_{\L_N}\x)\rangle)\le w_d\frac{2}{d+2}\bigl(\sqrt{d}\bigr)^d\E
^2(\xi(0))
\]
and
%
%
\begin{eqnarray}
\label{1211}
\frac{1}{|\L_N|}\mathop{\sum_{x,y\in\L_N}}_{|x-y|=1}\bigl((x-y)\cdot
u\bigr)^2&=&2|u|^2   \quad \mbox{and}
\nonumber
\\[-8pt]
\\[-8pt]
\nonumber
\frac{1}{|\L_N|}\mathop{\sum_{x\in\L,y\in\partial\L_N}}_{|x-y|=1}\bigl((x-y)\cdot u\bigr)^2&\le&
\frac{|u|^2}{N}\rightarrow0\qquad \mbox{as } N\rightarrow\infty
\end{eqnarray}
and that by standard SLLN arguments for i.i.d. random variables with
finite second moments
%
%
\begin{eqnarray}
\label{a11}
\frac{\sum_{x\in\L_N}\xi(x) (x\cdot u)}{N^{d+2}}\le|u|\frac{\sum
_{x\in
\L_N}|\xi(x)|}{N^{d+1}}\rightarrow0
\nonumber
\\[-8pt]
\\[-8pt]
\eqntext{\mbox{a.s. and in } L^1 \mbox
{ as } N\rightarrow\infty.}
\end{eqnarray}
\item[(b)] The statement follows from (\ref{lowerbound}), (\ref
{upbound}), (\ref{1211}) and Proposition~\ref{sllngree} by noting that
for very large $N$
\[
\hspace*{90pt}\frac{1}{N^d}\E(\langle\x,(G_{\L_N}\x)\rangle)=\frac
{2w_d}{d-2}\E(\xi^2(0)).\hspace*{90pt}\qed
\]
\end{longlist}
\noqed\end{pf*}

\begin{rem}
\label{surf12}
Note that due to the properties of the Green's function, for $d=1,2$ we
have that $\langle\xi, G_{\L_N} \xi\rangle_{\L_N}/|\L_N|$
diverges as
$N\rightarrow\infty$, and therefore, by the same reasoning as in
Theorem~\ref{non-exist} above, the surface tension does not exist for $d=1,2$.
\end{rem}

\subsubsection{\texorpdfstring{Existence of the surface tension when $\E(\xi(0))=0$}
{Existence of the surface tension when E(xi(0)) = 0}}
\label{tsf2}

In this section we prove Theorem~\ref{subadit}.
We start with a lemma which allows us to integrate out one height
variable $\phi(x)$ conditional upon the heights of its nearest neighbors.
\begin{lem}
\label{boundv}
Let the function $V$ satisfy (\ref{tag2}) and (\ref{tag22}). Then there
exists some constant $C>0$ such that for all $\gamma\in\R$, and all
$\phi(x),\xi(x)\in\RR$, $x\in\zd$, we have
%
%
\begin{eqnarray}
\label{boundveq}
&&\int_{\R}\exp\biggl[-\frac{1}{2} \sum_{y\in\zd, |y-x|=1}V\bigl(\phi
(y)-\phi
(x)\bigr)+\xi(x)\phi(x)\biggr]\,\ormd\phi(x)
\nonumber
\\[-8pt]
\\[-8pt]
\nonumber
&&\qquad\ge C\exp\biggl[-\frac{1}{2} \sum
_{y\in\zd,|y-x|=1}V\bigl(\phi(y)-\gamma\bigr)+\xi(x)\gamma\biggr].
\end{eqnarray}
\end{lem}

The proof of Lemma~\ref{boundv} closely follows the proof of
Lemma II.1 in~\cite{FS} and will be omitted.

Recall from (\ref{surftensdef}) that for any $\Lambda\in\Z^d$ and for
any fixed $u\in\R^d$
\[
Z^{{{{\psi}}_u}}_{\L}[\xi]=\int_{\R^\Lambda} \exp(-H_\Lambda
^{\psi
_u}[\xi])\,\ormd\phi_{\Lambda_N}.
\]
Let $a,l\in Z^d, a=(a_i)_{1\le i\le d}, l=(l_i)_{1\le i\le d}$ and let
$l_1'\in\Z,$ with $a_1<l_1'<l_1$.
We are going to prove an approximate subadditive relation for $-\log
Z^{\psi_u}_\Lambda$, where $\Lambda$ is taken to be the rectangle
$\bar
\Lambda^{a,l}$, as defined in (\ref{lal}), which is divided into three
rectangles by restricting the first coordinate to $[a_1, l'_1-1]$, $\{
l'_1\}$, and $[l'_1+1, l]$, respectively (see Figure~\ref{fig1}). To simplify the notation, we
denote for any $a, l\in\Z^d$ and $u,v\in\Z$
\[
\bar\Lambda^{a,l}_{[u, v]} := \Lambda_{[u,v]\times[a_2, l_2]\times
\cdots\times[a_d, l_d]}\quad  \mbox{and} \quad \bar\Lambda^{a,l}_{u} :=
\Lambda
_{\{u\}\times[a_2, l_2]\times\cdots\times[a_d, l_d]}.
\]
%
Using the above decomposition, we will derive in Lemma \ref
{surftensineq} the following formula:
%

%
\begin{figure}
\includegraphics{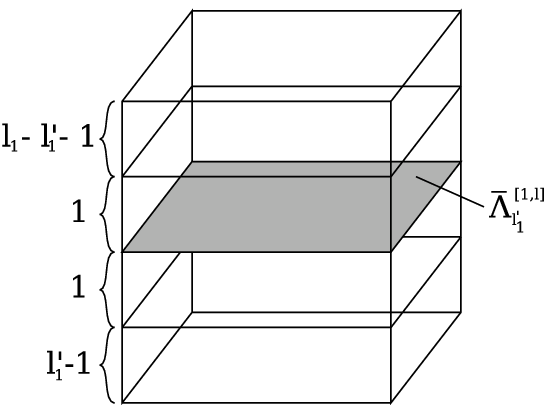}
\caption{}\label{fig1}
\vspace*{-3pt}
\end{figure}

\begin{lem}
\label{surftensineq}
Let the function $V$ satisfy (\ref{tag2}) and (\ref{tag22}). Then with
the notation above, we have for some $C>0$ and for $a_1\ge l_1+2$
%
%
\begin{eqnarray}
\label{subad1}
-\log(Z_{\bar\Lambda^{a,l}}^{{{\psi}}_u}[\xi])&\le&-\log
\bigl(Z_{ \bar\Lambda^{a,l}_{[a_1, l_1'-1]}}^{{{\psi}}_u}[\xi]\bigr)-\log
\bigl(Z_{\bar\Lambda^{a,l}_{[l_1'+1, l_1]}}^{{{\psi}}_u}[\xi]
\bigr)\nonumber\\[-3pt]
&&{}-\prod_{i=2}^d(l_i-a_i+1)\Biggl(\log C- \sum_{i=2}^d V(u_i)
\Biggr)\\[-3pt]
&&{}-\sum_{x\in\bar\Lambda^{a,l}_{l_1'}}u\cdot(l_1',x_2,\ldots
,x_d)\xi(x).\nonumber\vspace*{-3pt}
\end{eqnarray}
\end{lem}

\begin{pf}
We label the points $x\in\Lambda_{\bar\Lambda^{a,l}_{l_1'}}$ as odd or
even, depending on whether $\sum_{i=1}^d x_i$ is an odd or an even
number. We will bound $Z^{{{{\psi}}_u}}_{\bar\Lambda^{a,l}}[\xi]$
\textit{from below} by a product of $Z^{{{{\psi}}_u}}_{\bar\Lambda
^{a,l}_{[a_1, l_1'-1]}}[\xi]$, of $Z^{{{{\psi}}_u}}_{\bar\Lambda
^{a,l}_{[l_1'+1, l_1]}}[\xi]$ and of terms coming from\vspace*{-2pt} integrating out
the contribution of the elements of $\bar\Lambda^{a,l}_{l_1'}$ in
$H_{\bar\Lambda^{a,l}}^{\psi_u}[\xi](\phi).$\vspace*{1.5pt}
To do this, we will first integrate out the height variables at the odd
points in $\bar\Lambda^{a,l}_{l_1'}$ from $Z^{{{{\psi}}_u}}_{\bar
\Lambda
^{a,l}}[\xi]$ and then the even ones. We will do this by means of
Lemma~\ref{boundv}\vspace*{1pt} and by splitting $H_{\bar\Lambda^{a,l}}^{\psi
_u}[\xi
](\phi)$ into sums of potentials $V(\phi(x)-\phi(y))$, depending on
whether $x$ and $y$ belong to $\bar\Lambda^{a,l}_{[a_1, l_1'-1]}$,
$\bar
\Lambda^{a,l}_{[l_1'+1, l_1]}$, $\bar\Lambda^{a,l}_{l_1'}$ or
$\partial
\bar\Lambda^{a,l}$. Then by Lemma~\ref{boundv}, for each height
variable $\phi(x)$, $x\in\bar\Lambda^{a,l}_{l_1'}$ with $x$ odd,
(\ref{boundveq}) holds with $\gamma=u\cdot(l_1',x_2,\ldots,x_d)$ (we recall
that the boundary conditions for the two subdomains have the same tilt
$u$ as for the original domain). Explicitly, for each height variable
$\phi(x)$, $x\in\Lambda_{\bar\Lambda^{a,l}_{l_1'}}$ with $x$ odd,
we have
%
%
\begin{eqnarray}
\label{surf1}
&&\int_\R\exp\biggl[-\frac{1}{2}\sum_{j\in I}V\bigl(\phi(x+e_j)-\phi(x)\bigr)+\xi
(x)\phi(x)\biggr]\,\ormd\phi(x)
\nonumber
\\[-10pt]
\\[-10pt]
\nonumber
&&\qquad\ge C\exp\biggl[-\frac{1}{2}\sum_{j\in
I}V\bigl(\phi(x+e_j)-x\cdot u\bigr)+\xi(x) (x\cdot u)\biggr],
\end{eqnarray}
where $I:= \{\pm1,\pm2,\ldots,\pm d\}$. The point here is that Lemma
\ref{boundv} allows us to replace a height variable $\phi(x)$ by a
deterministic value $\gamma$.
Next we repeat the same procedure for each height variable $\phi(x)$,
$x\in\Lambda_{\bar\Lambda^{a,l}_{l_1'}}$ and $x$ even; since all
$\phi
(x+e_j)$, with $x+e_j\in\Lambda_{\bar\Lambda^{a,l}_{l_1'}}$ odd
nearest neighbors of $x$, have already been integrated out by (\ref
{surf1}), we have
%
\begin{eqnarray}
\label{surf2}
&&\int_\R\exp\biggl[-\frac{1}{2}\mathop{\sum_{j\in I}}_{j\neq\pm
1}V\bigl((x+e_j)\cdot u-\phi(x)\bigr)-V\bigl(\phi(x+e_1)-\phi(x)\bigr)\nonumber\\[-3pt]
&&\hspace*{116pt}{}-V\bigl(\phi
(x+e_{-1})-\phi
(x)\bigr)+\xi(x)\phi(x)\biggr]\,\ormd\phi(x)
\nonumber
\\[-10pt]
\\[-10pt]
\nonumber
&&\qquad\ge C\exp\Biggl[-\frac{1}{2}V\bigl(\phi(x+e_1)-x\cdot u\bigr)-\frac{1}{2}V\bigl(\phi
(x-e_1)-x\cdot u\bigr)\\[-3pt]
&&\hspace*{177pt}{}-\sum_{i=2}^d V(u_i)+\xi(x) (x\cdot
u)\Biggr].\nonumber
\end{eqnarray}
From (\ref{surf1}) and (\ref{surf2}) we get
\begin{eqnarray*}
Z^{{{{\psi}}_u}}_{\bar\Lambda^{a,l}}[\xi]&\ge& Z^{{{{\psi
}}_u}}_{\bar
\Lambda^{a,l}_{[a_1, l_1'-1]}}[\xi]  Z^{{{{\psi}}_u}}_{\bar\Lambda
^{a,l}_{[l_1'+1, l_1]}}[\xi]\\[-3pt]
&&{}\times \exp\Biggl(|\bar\Lambda^{a,l}_{l_1'}| \log
C-|\bar\Lambda^{a,l}_{l_1'}|\sum_{i=2}^d V(u_i)+\sum_{x\in{\bar
\Lambda
^{a,l}_{l_1'}}}\xi(x) (x\cdot u)\Biggr).
\end{eqnarray*}
%
Plugging $|\bar\Lambda^{a,l}_{l_1'}|=\prod_{i=2}^d(l_i-a_i+1)$ in the
above, we get (\ref{subad1}).
\end{pf}

\begin{pf*}{Proof of Theorem \protect\ref{subadit}} We will use Lemma
\ref{surftensineq} together with Proposition~\ref{akk} to prove in part
(a1) below that $\lim_{ N\rightarrow\infty}\sigma_{\Lambda
_{N}}[\xi
](u)$ exists for $\P$-almost all~$\xi$ and Lemma~\ref{surftensineq} and
Proposition~\ref{schurl1} to derive in part (a2) the $L^1$ convergence.
We will then use the a.s. and $L^1$ convergence in order to show in
part (b) that the surface tension is independent of the disorder $(\xi
(x))_{x\in\Z^d}$.

\begin{longlist}[(a1)]
\item[(a1)] We first need to rewrite (\ref{subad1}) in Lemma \ref
{surftensineq} in a form such that we can apply Proposition~\ref{akk}.
Let $a,l\in\Z^d, a=(a_i)_{1\le i\le d},l=(l_i)_{1\le i\le d}$, with
$a_i<l_i$ for $1\le i\le d$, be arbitrary and let, with the notation
from Lemma~\ref{surftensineq},
\[
g_{\bar\Lambda^{a,l}} :=\prod_{i=1}^d(l_i-a_i+1)\Biggl(\log C- \sum
_{i=1}^d \frac{V(u_i)}{2}\Biggr).\vadjust{\goodbreak}
\]
Let $l+1=(l_i+1)_{1\leq i\leq d}$ and define $\bar\Lambda^{a, l+1}$ as
in (\ref{lal}). Let
\[
f_{\bar\Lambda^{a, l+1}}[\xi](u):=-\log(Z_{\bar\Lambda^{a,
l}}^{\psi_u}[\xi])+\sum_{x\in\bar\Lambda^{a, l}}(u\cdot x) \xi
(x)+g_{\bar\Lambda^{a, l}}.
\]
Then from (\ref{subad1}) we have the following subadditivity formula
for $l_1\ge a_1+2$:
%
%
\begin{equation}
\label{subad2}
f_{\bar\Lambda^{a, l+1}}[\xi](u)\le f_{\bar\Lambda^{a,l+1}_{[a_1,
l_1']}}[\xi](u)+f_{\bar\Lambda^{a,l+1}_{[l_1'+1,l_1+1]}}[\xi](u).
\end{equation}
To get the subadditivity formula (\ref{subad2}) for all $l_1>a_1$, we
use an argument similar to the one we used to obtain (\ref{upbound}),
to bound for $l_1\in\{a_1,a_1+1\}$:
\begin{eqnarray*}
-\log Z_{\bar\Lambda^{a,l}_{[a_1,l_1]}}^{\psi_u}[\xi]&\le& \prod
_{i=2}^d(l_i-a_i+1)\bigl(2d V(0)+\sigma^{C_2/2}[\xi=0](u=0)\bigr)\\
&&{}-\sum
_{x\in\bar\Lambda^{a,l}_{[a_1,l_1]}}(u\cdot x) \xi(x)-\bigl\langle\xi, G_{\bar\Lambda^{a,l}_{[a_1,l_1]}} \xi\bigr\rangle_{\bar
\Lambda^{a,l}_{[a_1,l_1]}},
\end{eqnarray*}
%
where $\s^{C_2/2}[\xi=0](u=0)$ is defined as in Theorem \ref
{non-exist}(b). Taking into account that for all $\Lambda\in\Z^d$,
$\langle\xi, G_{\Lambda} \xi\rangle_{\Lambda}\ge0$, and making the
convention that for all $a_1\in\Z$
\begin{eqnarray*}
f_{\Lambda_{\bar\Lambda^{a,l}_{[a_1,a_1+1]}}}[\xi](u):&=&\prod
_{i=2}^2(l_i-a_i+1)\bigl(2d V(0)+|\sigma^{C_2/2}[\xi=0](u=0)
|\bigr)\\
&&{}-\sum_{x\in\bar\Lambda^{a,l}_{[a_1,a_1+1]}}(u\cdot x) \xi(x)+\bigl\langle\xi, G_{\bar\Lambda^{a,l}_{[a_1,a_1+1]}} \xi\bigr\rangle
_{\bar
\Lambda^{a,l}_{[a_1,a_1+1]}},
\end{eqnarray*}
it follows that for all $l_i>a_i,i=1,2,\ldots, d$, $f_{\bar\Lambda^{a,
l+1}}[\xi](u)$ satisfies the subadditivity property (\ref{subad2}) as
defined in
Proposition~\ref{akk}(b). We will check next that $f_{\bar\Lambda^{a,
l+1}}[\xi](u)$ satisfies conditions (a) and (c) of Proposition~\ref{akk}.
Recall that for $z\in\zd$, $\tau_{z}\phi(x)=\phi(x-z) \mbox
{for} x\in\zd
 \mbox{and} \phi\in\RR^{\zd}$. As $(\xi(x))_{x\in\zd}$ are
i.i.d., it
is easy to see that condition (a) of Proposition~\ref{akk} is
satisfied. We will
show next that (c) from Proposition~\ref{akk} also holds. Using the
lower bound in (\ref{1a}) and the fact that $\E(\xi(0))=0$, we have
that $f_{\bar\Lambda^{a, l+1}}[\xi](u)\in L^1$. Moreover, by the same
reasoning as that used to get (\ref{1a}), we have
\[
\frac{\E(f_{\bar\Lambda^{a, l+1} }[\xi](u))}{|\bar\Lambda
^{a, l}|}>\s^A[\xi=0](u=0)
-\frac{\E(\xi^2(0))\sum_{x\in\bar\Lambda^{a, l}}G_{\bar\Lambda^{a,
l}}(x,x)}{\bar\Lambda^{a,l}}-2d B.
\]
Since by Proposition~\ref{properties} we have that $\lim_{\Lambda\in
\Z
^d,|\Lambda|\uparrow\infty}G_{\Lambda}(x,x)=G(0,0)<\infty$, it
follows that
%
%
\begin{equation}
\label{bound}
\mathop{\inf_{a,l\in\Z^d,a_i<l_i}}_{ i=1,\ldots,d}\frac{\E(f_{\bar
\Lambda
^{a, l}}[\xi](u))}{|\bar\Lambda^{a, l}|}>-\infty
\end{equation}
and thus condition (c) of Proposition~\ref{akk} is also satisfied. It
follows that
%
%
\begin{equation}
\label{convf}
\lim_{N\rightarrow\infty} \frac{f_{\Lambda_{N}[\xi
](u)}}{N^d} \mbox{ exists a.s.}
\end{equation}
Together with (\ref{a11}) this proves that $\lim_{ N\rightarrow
\infty
}\sigma_{\Lambda_{N}}[\xi](u)$ exists for $\P$-almost all $\xi$.

\item[(a2)] To prove that $\lim_{ N\rightarrow\infty}\sigma
_{\Lambda
_{N}}[\xi](u)$
exists in $L^1$,
we will show that $f_{\bar\Lambda^{a, l+1}}[\xi](u)$ satisfies the
assumptions of Proposition~\ref{schurl1}. Note first that assumption
(a) is automatically satisfied, due to the subadditivity property
derived in~(\ref{subad2}). Similarly, assumption (d) is satisfied
because of (\ref{bound}). We will next prove that (b), (c) and (e) from
Proposition~\ref{schurl1} also hold. Let $z\in\Z^d$ and denote by
$(\tilde\psi)_u^z(x):=\sum_{i=1}^d(x_iu_i+z_i)$ for $x\in\partial
(\bar\Lambda^{a, l}+z)$. Then
%
%
\begin{eqnarray}
\label{station}
\qquad f_{\bar\Lambda^{a, l+1}+z}[\xi](u)&=&-\log\bigl(Z_{\bar\Lambda^{a,
l}+z}^{(\tilde\psi)_u^z}[\xi]\bigr)+\sum_{x\in\bar\Lambda^{a,
l}+z}(u\cdot x) \xi(x)+g_{\bar\Lambda^{a, l}+z}\\
&=&-\log(Z_{\bar\Lambda^{a, l}+z}^{\psi_u}[\xi])+\sum_{x\in
\bar\Lambda^{a, l}}(u\cdot x) \xi(x+z)+g_{\bar\Lambda^{a, l}},
\end{eqnarray}
where in the first equality we made in the integral formula for
$Z_{\bar
\Lambda^{a, l}+z}^{(\tilde\psi)_u^z}[\xi]$ the change of variables
${\hat\phi}(x):=\phi(x)+\sum_{i=1}^dz_iu_i$ for all $x\in\bar
\Lambda
^{a, l}+z$, and we used $g_{\bar\Lambda^{a, l}+z}=g_{\bar\Lambda^{a,
l}}$. Since $(\xi(x))_{x\in\Z^d}$ are i.i.d., (\ref{station}) proves
that (b), (c) and (e) from Proposition~\ref{schurl1} hold. It follows
that all assumptions of Proposition~\ref{schurl1} are satisfied. Therefore
\[
\frac{f_{\Lambda_{N}}[\xi](u)}{N^d} \mbox{ converges in }  L^1.
\]
Together with (\ref{a11}) this proves that $\lim_{ N\rightarrow
\infty
}\sigma_{\Lambda_{N}}[\xi](u)$ exists in $L^1$.

\item[(b)] Since we were unable to find in the literature a result for
multiparameter subadditive processes which we can apply directly as in
(a1) and (a2) to show that $\sigma(u)[\xi]$ is independent of the
disorder $\xi$, we will briefly sketch next a proof of the statement
for our case. For simplicity of notation, we restrict ourselves to
proving~(b) for $\Lambda_{[0,N]}$, where we recall that $\Lambda
_{[0,N]}=[0,N]^d\cap\Z^d$.

Let $k,n,r\in\Z_+$ such that $r<n$ and such that $N=kn+r$. For
$a=(a_i)_{1\le i\le d}\in\Z^d$, let $I_{a,n}:=\Lambda
_{[(a_1-1)n,a_1n]\times
\cdots\times[(a_d-1)n,a_dn]}$ and let $J_{N,k,n}^s:=\{z\in\Z
^d\dvtx kn\le
z_s\le N, 0\le z_i\le N \mbox{ for } i=\{1,2,\ldots,d\}\setminus\{s\}\}
$, where $s=1,2,\ldots, d$. Then
\[
\Lambda_{[0,N]}=\bigcup_{\{1\le a_i\le k,i=1,\ldots d\}}I_{a,n}\cup
\bigcup_{\{1\le s\le d\}}J_{N,k,n}^{s}.
\]
In words, we are partitioning $\Lambda_{[0,N]}$ into the union of cubes
of side lengths~$n$, which are the $I$'s, and the $J$'s represent the
leftover boundary terms because~$N$ may not be divisible by $n$. Thus
written, $\Lambda_{[0,N]}$ is a union of disjoint sets. From repeated
application of (\ref{subad2}), we have
%
%
\begin{equation}
\label{repsubad2}
f_{\Lambda_{[0,N]}}[\xi](u)\le\sum_{\{1\le a_i\le k,i=1,\ldots, d\}
}f_{I_{a,n}}[\xi](u)+\sum_{s=1}^d f_{J_{N,k,n}^{s}}[\xi](u).
\end{equation}
The key of the proof is that we can use the ergodic theorem for the
first sum in the right-hand side in (\ref{repsubad2}) and that the
boundary terms coming from the $J$'s are negligible. Combining this
with the a.s. and the $L^1$ convergence of $N^{-d}f_{\Lambda
_{[0,N]}}[\xi](u)$ proved in (a1) and (a2), the proof follows now
similar steps to the proof of Theorem 1.10 from~\cite{liggett} and will
be omitted.\qquad\qed
\end{longlist}
\noqed\end{pf*}

\subsection{Existence of shift-covariant random gradient Gibbs measures
with given tilt}
\label{esvg1}

This subsection is structured as follows: in Section~\ref{esvg1} we
construct in~(\ref{heldfixed}) a sequence of spatially averaged
finite-volume gradient Gibbs measures $(\bar\mu^{u}_{\L}[\xi]
)_{\L\subset\Z^d}$, such that $(\int\P(\ormd\xi)\bar\mu
^{u}_{\L
}[\xi])_{\L\subset\Z^d}$ is tight, as shown in Proposition \ref
{averaged}, and shift-invariant. In Section~\ref{esvg2} we will use
the tightness of $(\int\P(\ormd\xi)\bar\mu^{u}_{\L}[\xi])_{\L
\subset\Z^d}$ to prove in Theorem~\ref{thm1} the existence of a
shift-covariant random gradient Gibbs measure with a given tilt $u\in
\R^d$.

\subsubsection{Tightness of the averaged measure}

In order to prove tightness of the finite-volume gradient Gibbs
measures averaged over the disorder, we look at the finite-volume Gibbs measures
with tilt $u\in\R^d$ and boundary condition $\psi_u(x)=u\cdot x$:
\begin{eqnarray}
\label{1eqn1}
\nu^{\psi_u}_{\L}[\xi](\ormd\phi)&=&\frac{1}{Z_{\L}^{\psi_u}[\xi
]}\exp\biggl(
-\frac{1}{2}\mathop{\sum_{x,y\in\L}}_{|x-y|=1}V\bigl(\phi(x)-\phi(y)\bigr)
\nonumber\\
&&\hspace*{56pt}{}-\mathop{\sum_{x\in\L,y\in\partial\L}}_{|x-y|=1}V\bigl(\phi(x)-\psi_u(y)\bigr)
\\
&&\hspace*{117pt}{}+\sum_{x\in\L}\xi(x) \phi(x)
\biggr)
\,\ormd\phi_{\L}\d_{\psi_u}(\ormd\phi_{\Z^d\setminus\L}).\nonumber
\end{eqnarray}
Let us look now at the quantity
%
%
\begin{eqnarray}
\label{1'}
\qquad F_{\b,u,\L}[\xi_{\L}]&:=& \log\int\nu^{\psi_u}_{\L}[\xi](\ormd
\phi)
\nonumber
\\[-8pt]
\\[-8pt]
\nonumber
&&{}\times\exp\biggl(
+\frac{\b}{2}\sum_{x,y\in\zd,|x-y|=1}\bigl(\phi(x)-\phi(y)-u\cdot(x-y)\bigr)^2
\biggr),
\end{eqnarray}
for $\b>0$ sufficiently small. In (\ref{1'}), the sum over $x,y\in
\zd
,|x-y|=1$, can be taken
to include all the bonds on $\Z^d$ due to the fact that $\phi=\psi_u$
on $\L^c$. Note that
$F_{\b,u,\L}$ is the difference between the original free energy in the
volume $\L$ and the free energy in the volume $\L$ where we have added
the term
$\frac{\b}{2}\sum_{x,y\in\zd,|x-y|=1}(\phi(x)-\phi(y)-u\cdot(x-y))^2$
to the Hamiltonian.

We first note the following disorder-dependent upper bound for $F_{\b
,u,\L}$.

\begin{lem}\label{Lemmaoneone}
Let $d\ge3$. Assume that $V$ satisfies (\ref{tag2}) and (\ref
{tag22}). Then
%
%
\begin{eqnarray}\label{1}
F_{\b,u,\L}[\xi_{\L}]&\leq&
-|\L|\bigl(\s_{\L}^{A-\b}[\xi=0](u=0) - \s_{\L}^{C_2/2}[\xi
=0](u=0)\bigr)\nonumber\\[-2pt]
&&{}+\mathop{\sum_{x,y\in\L\cup\partial\L}}_{|x-y|=1}
\bigl(B+V(0)\bigr)
\nonumber\\[-2pt]
&&{}-\frac{A-\b-C_2/2}{2}\mathop{\sum_{x,y\in\L\cup\partial\L}}_{|x-y|=1}\bigl((x-y)\cdot u\bigr)^2\\[-2pt]
&&{}+\frac{1}{2}\biggl(\frac{1}{A-\b}-\frac
{2}{C_2}\biggr)\langle\xi, G_{\L} \xi\rangle_{\L}\nonumber\\[-2pt]
&=:&\bar F_{\b,u,\L
}+\frac{\alpha }{2}\langle\xi, G_{\L} \xi\rangle_{\L},\nonumber
\end{eqnarray}
with the obvious definitions for $\bar F_{\b,u,\L}$ and $\alpha$.
\end{lem}

\begin{pf}
Using bounds $A s^2 -B\leq V(s)\leq V(0)+\frac{C_2}{2} s^2$ for the
potential~$V$, we have
%
%
\begin{eqnarray}\label{1}
&&\exp( F_{\b,u,\L}[\xi_{\L}])\nonumber\\[-2pt]
&&\qquad\leq\int\exp\biggl(
-\frac{1}{2}\mathop{\sum_{x,y\in\L}}_{|x-y|=1}\bigl(A\bigl(\phi(x)-\phi(y)\bigr)^2
-B\bigr)
\nonumber\\[-2pt]
&&\hspace*{64pt}{}-\mathop{\sum_{x\in\L,y\in\partial\L}}_{|x-y|=1}\bigl(A\bigl(\phi(x)-\psi
(y)\bigr)^2 -B\bigr)
+\sum_{x\in\L}\xi(x) \phi(x)
\biggr) \nonumber\\[-2pt]
&&\hspace*{8pt}\qquad\quad{} \times\exp\biggl(
+\frac{\b}{2}\sum_{x,y\in\zd,|x-y|=1}\bigl(\phi(x)-\phi(y)-u\cdot
(x-y)\bigr)^2
\biggr)\,\ormd\phi_{\L}\\[-2pt]
&&\qquad\quad\Big/
\int\exp\biggl(
-\frac{1}{2}\mathop{\sum_{x,y\in\L} }_{|x-y|=1}\biggl(\frac{C_2}{2}\bigl(\phi
(x)-\phi(y)\bigr)^2+V(0)\biggr)
\nonumber\\[-2pt]
&&\hspace*{78pt}{}-\mathop{\sum_{x\in\L,y\in\partial\L}}_{ \atop|x-y|=1}\biggl(\frac{C_2}{2}\bigl(\phi
(x)-\psi(y)\bigr)^2+V(0)\biggr)\biggr)\nonumber\\
&&\qquad\hspace*{22pt}\quad{} \times\exp\biggl(\sum_{x\in\L}\xi(x) \phi(x)
\biggr)\,\ormd\phi_{\L}.\nonumber
\end{eqnarray}
This, by the same reasoning as in the proof of Theorem~\ref{non-exist},
is equal to
%
%
\begin{eqnarray}\label{1}
\hspace*{-4pt}&&\int\exp\biggl(
-\frac{1}{2}\mathop{\sum_{x,y\in\L}}_{|x-y|=1}
\bigl((A-\b)\bigl(\tilde\phi(x)-\tilde\phi(y)\bigr)^2\nonumber\\[-3pt]
&&\qquad\hspace*{62pt}{} + (A-\b)\bigl((x-y)\cdot
u\bigr)^2-B\bigr)\nonumber\\[-3pt]
&&\hspace*{11pt}\qquad{} -\mathop{\sum_{x\in\L,y\in
\partial\L}}_{|x-y|=1}
\bigl((A-\b)(\tilde\phi(x))^2 + (A-\b)\bigl((x-y)\cdot u\bigr)^2-B\bigr)\nonumber\\[-3pt]
&&\hspace*{198pt}\qquad{}+\sum
_{x\in\L}\xi(x) \tilde\phi(x)
\biggr) \,\ormd\tilde\phi_{\L}
\nonumber
\\[-9pt]
\\[-9pt]
\nonumber
&&\qquad{}\Big/ \int\exp\biggl(
-\frac{1}{2}\mathop{\sum_{x,y\in\L}}_{|x-y|=1}
\biggl(\frac{C_2}{2}\bigl(\tilde\phi(x)-\tilde\phi(y)\bigr)^2\\[-3pt]
&&\hspace*{100pt}\qquad{}+ \frac
{C_2}{2}\bigl((x-y)\cdot u\bigr)^2+V(0)\biggr)\nonumber\\[-3pt]
&&\hspace*{42pt}\qquad{} -\mathop{\sum_{x\in\L
,y\in\partial\L} }_{|x-y|=1}
\biggl(\frac{C_2}{2}(\tilde\phi(x))^2 + \frac{C_2}{2}\bigl((x-y)\cdot
u\bigr)^2+V(0)\biggr)\nonumber\\[-3pt]
&&\hspace*{207pt}\qquad{}+\sum_{x\in\L}\xi(x) \tilde\phi(x)
\biggr) \,\ormd\tilde\phi_{\L},\nonumber
\end{eqnarray}
where we note the cancellation of a sum over $\x$'s and where, as in
the proof of Theorem~\ref{non-exist}, for all $x\in\L$ we used the
change of variables $\phi(x)=\tilde\phi(x)+x\cdot u$.
The statement of the Lemma follows now by computing the Gaussian
integrals above as in the proof of Theorem~\ref{non-exist}.
\end{pf}

Take $\rho_u(b):=\nabla\psi_u(b)$ for all $b\in(\zd)^*$ and consider
the corresponding gradient Gibbs measure $\mu_\L^{\rho_u}[\xi]$ as
given by (\ref{fingradgibbs}). Let us now define the \textit{spatially
averaged} measure $\bar\mu^{u}_{\L}[\xi]$ on gradient configurations
obtained by
%
%
\begin{equation}
\label{heldfixed}
\bar\mu^{u}_{\L}[\xi]:=\frac{1}{|\L|}\sum_{x\in\L}\mu^{\rho
_u}_{\L
+x}[\xi],
\end{equation}
where we recall that $\L+x:=\{z+x\dvtx z\in\L\}$. This is an extension to
our disorder-dependent case
of the construction on Gibbs measures with symmetries given in~\cite
{giorgii}, in formula (5.20) from Chapter $5.2$; the construction in~\cite{giorgii}
was used\vadjust{\goodbreak}
to get shift-invariant Gibbs measures. We note that in (\ref
{heldfixed}), the random field variables $\x$ are held fixed while the
volumes $\L+x$
are shifted around. We will first use the fact that the measure $
(\int\P(\ormd\xi)\bar\mu^{u}_{\L}[\xi])(\ormd\phi)$ is
shift-invariant in the proof of Proposition~\ref{averaged} below.
Then we will use $\bar\mu^{u}_{\L}[\xi]$ to construct shift-covariant
gradient Gibbs measures in Section~\ref{esvg2} by performing a further
average over the volumes.

In preparation for the proof of existence of random shift-covariant
gradient Gibbs measures, we will prove the following result on the
tightness of the family of averaged finite-volume random $\nabla\phi$-Gibbs
measures, and therefore on the existence, of the
infinite-volume $\nabla\phi$-Gibbs measures averaged over the
disorder.\vspace*{-3pt}

\begin{prop}
\label{averaged}
Suppose that $d\geq3$ and $\E(\xi(0))=0$. Assume that $V$ satisfies
(\ref{tag2}) and (\ref{tag22}). Then there exists a constant $K>0$ such
that for all $x_0, y_0\in\Z^d$ with $|x_0-y_0|=1$, the measure
\[
P^u_{\L}(\ormd\phi)
:=\biggl(\int\P(\ormd\xi)\bar\mu^{u}_{\L}[\xi]\biggr)(\ormd\phi)=
\biggl(\frac{1}{|\L|}\sum_{x\in\L}\int\P(\ormd\xi)\mu^{\rho_u}_{\L
+x}[\xi
]\biggr)(\ormd\phi)
\]
satisfies the estimate
%
%
\begin{equation}
\label{1}
\limsup_{N\uparrow\infty} P^u_{\L_N}\bigl[\bigl(\phi(x_0)-\phi(y_0)-u\cdot
(x_0 -y_0)\bigr)^2\bigr] \leq K.
\end{equation}
Hence the sequence of measures $P^u_{\L_N}$ is tight and thus possesses
a disorder-independent
limit measure (along subsequences of volumes) on gradient configurations.\vspace*{-3pt}
\end{prop}

\begin{pf} Let $f\dvtx \R^{\Z^d}\rightarrow[0,\infty)$ be given
by $f(\phi):=(\phi(x_0)-\phi(y_0)-u\cdot(x_0-y_0) )^2$; using
translation invariance of the
distribution of the disorder $(\xi(x))_{x\in\Z^d}$, we have
\begin{eqnarray*}
P^{u}_{\L}(f)&=&\biggl[\frac{1}{|\L|}\sum_{x\in\L}\E\mu^{\rho_u}_{\L
+x}[\xi]\biggr](f)=\frac{1}{|\L|}\sum_{x\in\L}(\E\mu^{\rho_u}_{\L
+x}[\xi])(f\circ\tau_x)
\nonumber
\\[-8pt]
\\[-8pt]
\nonumber
&=&\frac{1}{|\L|}\E\mu^{\rho_u}_{\L}[\xi]\biggl(\sum_{x\in\L}f\circ
\tau
_{x}\biggr).
\end{eqnarray*}
By the nonnegativity of $f$ we have for $\P$-almost all $\xi$
\begin{eqnarray*}
\label{1}
\mu^{\rho_u}_{\L}[\xi]\biggl(\sum_{x\in\L}f\circ\tau_{x}\biggr)&\leq&\mu
^{\rho_u}_{\L}[\x]
\biggl(\sum_{x,y\in\zd,|x-y|=1}\bigl(\phi(x)-\phi(y)-u\cdot(x-y)\bigr)^2
\biggr)\\[-3pt]
&=:&g[\xi].
\end{eqnarray*}
By writing $g[\xi]=({2}/{\beta}) \log e^{({\beta}/{2} )g[\xi
]}$ and
applying Jensen's inequality, we have
%
\begin{eqnarray*}
P^{u}_{\L}(f)
&\leq& \frac{1}{|\L|} \E\mu^{\rho_u}_{\L}[\x]
\biggl(\sum_{x,y\in\zd,|x-y|=1}\bigl(\phi(x)-\phi(y)-u\cdot(x-y)\bigr)^2\biggr)\\
&\leq& \frac{2}{\b|\L|} \E\log\mu^{\rho_u}_{\L}[\xi]\biggl(\exp
\biggl(\frac{\b}{2}\sum_{x,y\in\zd,|x-y|=1}\bigl(\phi(x)-\phi(y)-u\cdot
(x-y)\bigr)^2\biggr)\biggr).\vadjust{\goodbreak}
\end{eqnarray*}
By Lemma~\ref{Lemmaoneone} we get when $\L=\L_N$ the upper bound
%
%
\begin{equation}
P^{u}_{\L_N}(f)\leq\frac{2}{\b|\L_N|}\bar F_{\b,u,\L_N}
+ \frac{2}{\b|\L_N|} \E\biggl(\frac{\alpha }{2}\langle\xi, G_{\L_N}
\xi
\rangle_{\L_N}\biggr),
\end{equation}
which is bounded uniformly in $\L_N$, as $\frac{\bar F_{\b,u,\L
_N}}{|\L
_N|}$ is uniformly bounded by Theorem~\ref{subadit} and by (\ref
{1211}), and $0\le\frac{1}{|\L_N|} \E(\langle\xi, G_{\L_N} \xi
\rangle_{\L_N})\le G(0,0)+1$, by Proposition~\ref{properties}(ii)
and $\E(\xi(0))=0$. This proves the claim.
\end{pf}

\subsubsection{Existence of shift-covariant random gradient Gibbs
measures with given tilt}
\label{esvg2}

In this subsection we will prove our main result, Theorem~\ref{thm1},
of existence of a shift-covariant random gradient Gibbs measure $\hat
\mu
^u[\xi]$ with a given tilt $u\in\R^d$.
In the proof, we will first construct a candidate $\hat\mu^u[\xi]$ by
taking suitable subsequential weak limits, and then in two subsequent
Lemmas~\ref{bargradgibbs} and~\ref{shiftlemma}, we will prove,
respectively, that $\P$-a.s., our candidate~$\hat\mu^u[\xi]$ is a
gradient Gibbs measure, and is translation-covariant.

To construct a candidate $\hat\mu^u[\xi]$, we will need to perform a
further average of $\bar\mu^u[\xi]$ over the volumes $\Lambda$, and to
find a deterministic sequence $(m_r)_{r\in\N}$, along which there is a
weak limit for $\P$-a.e. $\xi$. This will be facilitated by Theorem 1a
from~\cite{KOM}, which we state below.
%
\begin{prop}
\label{kom}
If $(\zeta_n)_{n\in\N}$ is a sequence of real-valued random variables
with $\lim\inf_{n\rightarrow\infty}\E(|\zeta_n|)<\infty$, then there
exists a subsequence $\{\theta_n\}_{n\in\N}$ of the sequence $\{
\zeta
_n\}_{n\in\N}$ and an integrable random variable $\theta$ such that for
any arbitrary subsequence $\{\tilde{\theta}_n\}_{n\in\N}$ of the
sequence $\{\theta_n\}$, we have
\[
\lim_{n\rightarrow\infty}\frac{\tilde{\theta}_1+\tilde{\theta
}_2+\cdots
+\tilde{\theta}_n}{n}=\theta \qquad \P\mbox{-almost surely}.
\]
\end{prop}

We are now ready to prove the existence of shift-covariant gradient
Gibbs measures in Theorem~\ref{thm1}, which follows immediately from
the next Proposition.
\begin{prop}
\label{existgibbs}
Suppose that $d\geq3$ and $\E(\xi(0))=0$. Assume that $V$ satisfies
(\ref{tag2}) and (\ref{tag22}). Then there is a deterministic sequence
$(m_r)_{r\in\N}$ in~$\N$ such that for $\P$-almost every $\xi$,
%
%
\begin{equation}
\hat\mu^{u}_{k}[\xi]:=\frac{1}{k}\sum_{i=1}^k {\bar\mu}^{u}_{\L
_{m_{i}}}[\xi]
\end{equation}
converges as $k\to\infty$ weakly to $\hat\mu^u[\xi]$, which is a
shift-covariant random gradient Gibbs measure defined as in Definition
\ref{shiftcov1}.
\end{prop}

\begin{pf} We will prove first that there exists a
deterministic sequence $(m_r)_{r\in\N}$ in $\N$ such that $(\hat\mu
^{u}_{k}[\xi])_{k\in\N}$ converges a.s.\ to a random measure $\hat
\mu^{u}[\xi]$. We will then show that $\hat\mu^{u}[\xi]$ is a.s. a
gradient Gibbs measure, is translation-covariant
and that $\xi\rightarrow\hat\mu^{u}[\xi]$ is a measurable map.

Let $(f_i)_{i\in\N}$ be a countable collection of functions in
$C_b(\chi
)$, such that a~sequence of probability measures $\mu_n \in P(\chi)$
converges weakly to $\mu\in P(\chi)$ if and only if $\mu_n(f_i)\to
\mu
(f_i)$ for all $i\in\N$. Such a countable family $(f_i)_{i\in\N}$ in
$C_b(\chi)$ is explicitly given, for example, in the general setting of
separable and complete metric spaces in Proposition 3.17 from \cite
{RES} or in Lemma 1.1 from~\cite{Kal}. To show
that for a given sequence $(m_r)_{r\in\N}$ and a random measure $\hat
\mu
[\xi]$, $\hat\mu_k[\xi]$ converges a.s. to $\hat\mu[\xi]$, it suffices
to show that $\hat\mu_{k}[\xi](f_i) \to\hat\mu[\xi](f_i)$ almost
surely for each $i\in\N$.


For each $N\in\N$ and $x,y\in\Z^d$ with $|x-y|=1$, define
%
%
\begin{equation}
\label{heldfixed2}
X_{N;x,y}[\xi]:=\bar\mu_{\L_N}^u[\xi]\bigl(\bigl(\phi(x)-\phi(y)-u\cdot(x
-y)\bigr)^2\bigr).
\end{equation}
Take now the countable sequence containing both the family $(\bar
\mu_{\L_N}^u[\xi](f_i))_{i,N\in\N}$ and $(X_{N;x,y}[\xi])_{N\in
\N
, x,y\in\Z^d\atop|x-y|=1}$. We note that since $(f_i)_{i\in\N}$ are
bounded functions, $\lim\inf_{N\uparrow\infty}\E(\bar\mu_{\L
_N}^u[\xi](|f_i|))\hspace*{-0.3pt}<\hspace*{-0.3pt}|f_i|_{\infty}\hspace*{-0.3pt}<\hspace*{-0.3pt}\infty$.
Note also that $\liminf_N \E(X_{N;x,y}[\xi])\hspace*{-0.3pt} <\infty$ by
Proposition~\ref{averaged}. Therefore by Proposition~\ref{kom}, for
each $x_0,y_0\in\Z^d$ with $|x_0-y_0|=1$, there exists a sequence
$(n_r)_{r\in\N}$ and a random variable $\kappa_{x_0, y_0}$, both
depending on $x_0$ and $y_0$, such that
\[
\lim_{k\uparrow\infty}\frac{1}{k}\sum_{r=1}^k X_{n_r;x_0,y_0}[\xi
]=\kappa_{x_0,y_0}[\xi]\qquad  \mbox{for } \P\mbox{-almost every } \xi.
\]
Moreover
\[
\lim_{k\uparrow\infty}\frac{1}{k}\sum_{j=1}^k
X_{n_{r_j};x_0,y_0}[\xi
]=\kappa_{x_0,y_0}[\xi] \qquad \mbox{for } \P\mbox{-almost every } \xi
\]
holds also for every further subsequence $(n_{r_j})_{r_j\in\N}$ of
$(n_r)_{r\in\N}$. We take an arbitrary such subsequence $n_{r_j}$. By
Proposition~\ref{kom}, there exists a subsequence $(n'_r)_{r\in\N}$ of
$(n_{r_j})_{r_j\in\N}$ and a random variable $\rho_1$, both depending
on $x_0$ and $y_0$, such that
%
\[
\lim_{k\uparrow\infty}\frac{1}{k}\sum_{j=1}^k \bar\mu_{\L
_{n'_{r_j}}}^u[\xi](f_1) =\rho_{1}[\xi]\qquad  \mbox{for } \P\mbox{-almost
every } \xi.
\]
Moreover
\[
\lim_{k\uparrow\infty}\frac{1}{k}\sum_{j=1}^k \bar\mu_{\L
_{n''_{r_j}}}[\xi](f_1)=\rho_{1}[\xi]\qquad  \mbox{for } \P\mbox{-almost
every } \xi
\]
holds also for every further subsequence $n''_{r_j}$ of
$n'_{r_j}$.\vadjust{\goodbreak}

We repeat this procedure for each $x,y\in\Z^d,|x-y|=1$ and for each
\mbox{$i\in\N$}. By a Cantor diagonalization argument over the countably many
$x,y\in\Z^d,|x-y|=1$ and over the $i\in\N$, there exists a
deterministic sequence $(m_r)_{r\in\N}$ in $\N$ and random variables
$(\kappa_{x,y}[\xi])_{x,y\in\Z^d,|x-y|=1}$ and $(\rho_i[\xi
])_{i\in\N}$
such that for $\P$-almost every $\xi$,
%
%
\begin{eqnarray}
\label{limkom}
\lim_{k\uparrow\infty}\frac{1}{k}\sum_{r=1}^k X_{m_{r};x,y}[\xi
]&=&\kappa
_{x,y}[\xi] \quad \mbox{and}
\nonumber
\\[-8pt]
\\[-8pt]
\nonumber
  \lim_{k\uparrow\infty}\frac{1}{k}\sum
_{r=1}^k\bar\mu_{\L_{m_r}}^u[\xi](f_i)&=&\rho_i[\xi],
\end{eqnarray}
for all $x,y\in\Z^d$ and all $i\in\N$.
In particular, we get from (\ref{limkom}) that\break $\sup_{k\in\N}{\hat
X}_{k;x,y}[\xi]\le C(\kappa_{x,y}[\xi])$ for some $C(\kappa
_{x,y}[\xi
])>0$. Therefore for all $b\in(\Z^d)^*$, with $b=(x,y)$, we have for
$\P$-almost every $\xi$
\[
\lim_{L\uparrow\infty}\sup_{k\in\N}{\hat\mu}_k^u[\xi]\bigl(\eta
\dvtx |\eta(b)|\ge
L\bigr)\le\lim_{L\uparrow\infty}\sup_{k\in\N}\frac{{\hat
X}_{k;x,y}[\xi]}{L^2}=0.
\]
This means that for $\P$-a.s. all $\xi$, there exists a (possibly)
\textit{random} subsequence $(k'[\xi])$ such that $({\hat\mu
}_{k'[\xi]}^u[\xi])_{k'[\xi]}$ is tight and converges weakly to a
random measure $\hat\mu^u[\xi]$. The random subsequence $(k'[\xi])$ is
used only for tightness;
in fact the subsequence becomes nonrandom again as we return below to
the deterministic subsequence $(m_r)$. Moreover, we have ${\hat\mu
}_{k'[\xi]}^u[\xi](f_i)\rightarrow\hat\mu^u[\x](f_i)$ for all
$i\in\N$.
Due to (\ref{limkom}), and by the uniqueness of the limit point, we get
that $\rho_i[\xi]=\hat\mu^u[\x](f_i)$ for all $i\in I$. Since
$\hat\mu
_{k}^u[\xi](f_i)\rightarrow\hat\mu^u[\x](f_i)$, it follows that~$\hat
\mu_k^u[\xi]$ converges a.s. to a random measure $\hat\mu^u[\xi]$.

From Lemma~\ref{bargradgibbs} below, we get that for $\P$-almost all
$\xi$, $\hat\mu^u[\xi]$ is a gradient Gibbs measure and from Lemma
\ref
{shiftlemma} below, that $\hat\mu^u[\xi]$ is translation-covariant for
$\P$-almost all $\xi$.

It only remains to prove that $\xi\rightarrow\hat\mu^u[\xi]$ is a
measurable map. We recall that the disorder is defined on the
probability space $(\Omega, \mathcal{F}, \P)$. With a given tilted
boundary condition $\psi_u$, $\mu^{\psi_u}_{\Lambda}[\xi]$ is
clearly a
measurable function of the disorder field $\xi$. Since $\hat\mu^u$ is
constructed as a pointwise (w.r.t. $\xi$) limit of averages of such
measurable $\P(\chi)$-valued functions of $\xi$, $\hat\mu^u$ is
also a
measurable $\P(\chi)$-valued function of~$\xi$.
\end{pf}

We will prove next Lemmas~\ref{bargradgibbs} and~\ref{shiftlemma}. The
setup is as before; that is,~$\hat\mu^u_k[\xi]$ is defined as in
(\ref{heldfixed}), and the assumption is that along a deterministic
subsequence $(m_i)_{i\in\N}$ in $\N$, we have weak convergence of
$\hat
\mu^u_k[\xi]$ to $\hat\mu^u[\xi]$ for $\P$-almost all $\xi$.

%
\begin{lem}
\label{bargradgibbs}
For $\P$-almost all $\xi$, the limit $\hat\mu^u[\xi]$
is a gradient Gibbs measure.
\end{lem}

\begin{pf} In order to show that $\hat\mu^{u}[\xi]$ is a
gradient Gibbs measure, we have to show that for each fixed $\xi$, for
all $F\in C_b(\chi)$ and for all $J\subset\Z^d$ we have
%
%
\begin{equation}
\label{gibbseqn}
\int\hat\mu^u[\xi](\ormd\tilde\rho)\int\mu_{J}^{\tilde\rho
}[\xi](\ormd
\eta)F(\eta)=\int\hat\mu^u[\xi](\ormd\eta)F(\eta).
\end{equation}
Using the compatibility of the kernels, namely
\[
\int\mu^{\rho_u}_{\L}[\xi](\ormd\tilde\rho)\mu^{\tilde\rho
}_{J}[\xi]
=\mu^{\rho_u}_{\L}[\xi]\qquad  \mbox{for }  J\subset\L\subset\Z^d,
\]
we have
%
%
\begin{eqnarray}
\label{gradeqn}
&&\int\bar\mu^{u}_{\L}[\xi](\ormd\tilde\rho)\mu^{\tilde\rho
}_{J}[\xi
]\nonumber\\
&&\qquad=\frac{1}{|\L|}\sum_{x\in\Lambda}\int\mu^{\rho_u}_{\L
+x}[\xi](\ormd
\tilde\rho)\mu^{\tilde\rho}_{J}[\xi]
\nonumber\\
&&\qquad=\frac{1}{|\L|}\biggl(\sum_{x\in\L\dvtx  J\sb\L+x }+ \sum_{x\in\L\dvtx
J\not
\sb\L+x }\biggr)
\int\mu^{\rho_u}_{\L+x}[\xi](\ormd\tilde\rho)\mu^{\tilde\rho
}_{J}[\xi
]\\
&&\qquad=\frac{1}{|\L|}\sum_{x\in\L\dvtx  J\sb\L+x }\mu^{\rho_u}_{\L
+x}[\xi
]+\frac{1}{|\L|}
\sum_{x\in\L\dvtx  J\not\sb\L+x }
\int\mu^{\rho_u}_{\L+x}[\xi](\ormd\tilde\rho)\mu^{\tilde\rho
}_{J}[\xi
]\nonumber\\
&&\qquad=\bar\mu^{u}_{\L}[\xi]+ \frac{1}{|\L|}
\sum_{x\in\L\dvtx  J\not\sb\L+x }
\biggl(\int\mu^{\rho_u}_{\L+x}[\xi](\ormd\tilde\rho)\mu^{\tilde\rho
}_{J}[\xi]-\mu^{\rho_u}_{\L+x}[\xi]\biggr).\nonumber
\end{eqnarray}
%
Fix $J\subset\Z^d$ and take $k\in\N$ large enough. Applying (\ref
{gradeqn}) to the subsequence $(\L_{m_i})_{1\le m_i\le k}$ and to an
arbitrary $F\in C_b(\chi)$, we have
%
%
\begin{equation}
\label{graden11}\quad
\hat\mu^{u}_{k}[\xi](\mu^{\tilde\rho}_J[\xi](F))=\frac
{1}{k}\sum
_{i=1}^k\bar\mu^{u}_{\L_{m_i}}[\xi](F)+ \frac{1}{k}\sum
_{i=1}^kR(\L
_{m_i},J,F)[\xi],
\end{equation}
where $|R(\L_{m_i},J,f)[\xi]|\leq\frac{C(f)}{|\L_{m_i}|}
\sum_{x\in\L_{m_i}\dvtx  J\not\sb\L_{m_i}+x }1$, for all $1\le i\le k$ and
for some constant $C(f)>0$.
In order to prove (\ref{gibbseqn}), we need to take $k\rightarrow
\infty
$ on both sides of (\ref{graden11}). To do that, we have to prove first
that for all $F\in C_b(\chi)$ and for all fixed $J\subset\Z^d$ we have
%
%
\begin{equation}
\label{limitgibbs}
\int\hat\mu[\xi](\ormd\tilde\rho)(\mu^{\tilde\rho}_{J}[\xi](F))
=\lim_{k\uparrow\infty}\hat\mu^{u}_{k}[\xi](\mu^{\tilde\rho
}_J[\xi](F)).
\end{equation}
To show (\ref{limitgibbs}), it is sufficient to show that for all
$F\in
C_b(\chi)$ the function $\mu^{\tilde\rho}_{J}[\xi](F)\in C_b(\chi
)$ as
a function in $\tilde\rho$; then (\ref{limitgibbs}) will follow by the
hypothesis. The boundedness of $\mu^{\tilde\rho}_{J}[\xi](F)$ follows
immediately due to the boundedness of $F$. To prove continuity of $\mu
^{\tilde\rho}_{J}[\xi](F)$, fix $\tilde\rho\in\chi$ arbitrarily. As
$\chi$ equipped with the metric $\Vert_r$ is a complete metric space, we
can take now a~sequence $({\tilde\rho}_n)_{n\in\N}\in\chi$ such that
$\lim_{n\uparrow\infty}{\tilde\rho}_n=\tilde\rho$ in $\chi$; we
have to
show that $\lim_{n\uparrow\infty}\mu^{{\tilde\rho}_n}_{J}[\xi
](F)=\mu
^{\tilde\rho}_{J}[\xi](F)$. In view of the fact that $V\in C^2(\R
)$, we
note now that both the integrand in the numerator, and the integrand in
the denominator, of $\lim_{n\uparrow\infty}\mu^{{\tilde\rho
}_n}_{J}[\xi
](F)$ converge as ${\tilde\rho}_n\rightarrow{\tilde\rho}$; moreover,
due to the bounds $A s^2 -B\leq V(s)\leq V(0)+\frac{C_2}{2} s^2$ on the
potential $V$ and by a similar reasoning as in the proof of Lemma \ref
{Lemmaoneone}, these integrands are uniformly bounded by integrable
functions. Applying now Lebesgue's dominated convergence theorem
separately to the numerator and to the denominator gives $\lim
_{n\uparrow\infty}\mu^{{\tilde\rho}_n}_{J}[\xi](F)=\mu^{\tilde
\rho
}_{J}[\xi](F)$, and therefore (\ref{limitgibbs}) holds.
Taking $k$ to infinity in (\ref{graden11}) and using (\ref
{limitgibbs}), we get
%
\begin{eqnarray*}
\int\hat\mu[\xi](\ormd\tilde\rho)(\mu^{\tilde\rho}_{J}[\xi
](F))&=&\lim
_{k\uparrow\infty}\frac{1}{k}\sum_{i=1}^k\bar\mu^{u}_{\L
_{m_i}}[\xi
](F)+ \lim_{k\uparrow\infty}\frac{1}{k}\sum_{i=1}^k R(\L
_{m_i},J,F)[\xi
]\\
 &=&\hat\mu^{u}[\xi](F)+0,
\end{eqnarray*}
where the convergence holds due to the fact that $F\in C_b(R^{(\zd
)^*})$ and\break $\sum_{i=1}^k R(\L_{m_i}, J,F)[\xi]/k $ goes to zero
uniformly in $\x$, due to the upper bound on $|R(\L,J,F)[\xi]|$.
This proves that (\ref{gibbseqn}) holds.
\end{pf}

\begin{lem}
\label{shiftlemma}
For $\P$-almost all $\xi$, the limit $\hat\mu^u[\xi]$ is
translation-covariant, that is, for all $v\in\Z^d$ and for all $F\in
C_b(\chi)$, we have
%
%
\begin{equation}
\hat\mu^{u}[\xi](F\circ\tau_v)=\hat\mu^{u}[\tau_v\xi](F),
\end{equation}
where we recall that $(\tau_v\xi)(z)=\xi(z-v)$ for all $z\in\Z^d$.
\end{lem}

\begin{pf}
Fix $v\in\Z^d$. Then we have
%
%
\begin{eqnarray}
\label{siebenunddreissig1}
&&\hat\mu^{u}[\xi](F\circ\tau_{v})-\hat\mu^{u}[\tau_v\xi](F)
\nonumber\\
&&\qquad=\lim_{k\uparrow\infty}\frac{1}{k}\sum_{i=1}^k
\frac{1}{|\L_{m_i}|}\biggl(\sum_{x\in\L_{m_i}}\mu^{\rho_u}_{\L
_{m_i}+x}[\xi](F\circ\tau_{v})
\\
&&\hspace*{88pt}\qquad{}-\sum_{x\in\L_{m_i}}\mu^{\rho_u}_{\L_{m_i}+x}[\tau_v\xi](F)\biggr).\nonumber
\end{eqnarray}
The terms inside the last bracket equal
\begin{eqnarray*}
&&\sum_{x\in\L_{m_i}}\mu^{\rho_u}_{\L_{m_i}+x}[\xi](F\circ\tau_{v})
-\sum_{x\in\L_{m_i}}\mu^{\rho_u}_{\L_{m_i}+x}[\tau_v\xi
](F)\\
&&\qquad=\sum_{x\in
\L_{m_i+v}}\mu^{\rho_u}_{\L_{m_i}+x}[\xi](F)
-\sum_{x\in\L_{m_i}}\mu^{\rho_u}_{\L_{m_i}+x}[\xi](F).
\end{eqnarray*}
Most terms on the right-hand side cancel. Therefore, for a bounded
function~$F$ such that $\Vert F\Vert_{\infty} \leq C(F)$ for some
$C(F)>0$, we have
%
%
\begin{equation}
\label{2}
|\hat\mu^{u}[\xi](F)-\hat\mu^{u}[\tau_v\xi](F)|\leq
\lim_{k\uparrow\infty}\frac{C(F)}{k}\sum_{i=1}^k\frac{|{\L}_{m_i}
\triangle({\L}_{m_i} +v)|}{|{\L}_{m_i}|},
\end{equation}
where we denoted by $\Delta$ the symmetric difference of the sets $\L$
and $\L+v$.
But $|\L_{m_i} \triangle(\L_{m_i}+v)|$
goes to zero when divided by $|\L_{m_i}|$, uniformly in $m_i$, which
implies that (\ref{2})
goes to zero also.
This shows the translation-covariance.
\end{pf}

\begin{pf*}{Proof of Theorem \protect\ref{thm1}(\textup{a})} Proposition \ref
{existgibbs} implies the existence of a~random gradient Gibbs $\hat\mu
^u[\xi]$. We prove next that $\hat\mu^u[\xi]$ satisfies (\ref{tilt1}).
Given the tilt $u\in\R^d$, the limit $\hat\mu^u[\xi]$ we construct is
the weak limit of the~$\hat\mu_k^u[\xi]$. We next calculate what is the
expected tilt over a given bond under the measure $\hat\mu_k^u[\xi]$,
averaged over the disorder. For any $m_i$ in the deterministic sequence
$(m_i)_{1\le i\le k}$ and for $b_1:=(0,e_1)$, we have by means of~(\ref
{heldfixed}) and of Definition (\ref{fingradgibbs})
%
%
\begin{eqnarray}
\label{147}
\bar\mu_{m_i}^{\rho_u}[\xi](\eta(b_1))&=&\frac{1}{|\L_{m_i}|}
\sum
_{x\in\L_{m_i}}\mu^{\rho_u}_{\L_{m_i}+x}[\xi](\eta(b_1))\nonumber\\
&= &\frac
{1}{|\L
_{m_i}|} \sum_{x\in\L_{m_i}}\nu^{\psi_u}_{\L_{m_i}+x}[\xi]\bigl(\phi
(e_1)-\phi(0)\bigr)
\nonumber
\\[-8pt]
\\[-8pt]
\nonumber
&=&\frac{1}{|\L_{m_i}|} \sum_{x\in\L_{m_i}}\nu^{\tau_{-x}\psi
_u}_{\L
_{m_i}}[\tau_{-x}\xi]\bigl(\phi(e_1-x)-\phi(-x)\bigr)\\
&=&\frac{1}{|\L_{m_i}|} \sum_{x\in\L_{m_i}}\nu^{\psi_u}_{\L
_{m_i}}[\tau
_{-x}\xi]\bigl(\phi(e_1-x)-\phi(-x)\bigr),\nonumber
\end{eqnarray}
where for the third equality we made for all $y\in\L_{m_i}$ the change
of variables $\phi(y)\rightarrow\phi(y)+\sum_{i=1}^d u_i x_i$ under
each integral. Let
\begin{eqnarray*}
\bar\Lambda^{-m_i,m_i}_{m_i} &:=& \Lambda_{\{m_i\}\times[-m_i,
m_i]\times\cdots\times[-m_i, m_i]} \quad \mbox{and} \\
 \bar\Lambda
^{-m_i,m_i}_{-m_i} &:=& \Lambda_{\{-m_i\}\times[-m_i, m_i]\times\cdots
\times[-m_i, m_i]}.
\end{eqnarray*}
Averaging over the disorder in (\ref{147}), we get
\begin{eqnarray*}
&&\E\biggl(\int\bar\mu_{m_i}^{\rho_u}[\xi](\ormd\eta)\eta(b_1)
\biggr)\\
&&\qquad=\frac{1}{|\L_{m_i}|} \sum_{x\in\L_{m_i}} \E\biggl(\int\nu^{\psi
_u}_{\L_{m_i}}[\tau_{-x}\xi](\ormd\phi)\bigl(\phi(e_1-x)-\phi(-x)\bigr)\biggr)\\
&&\qquad=\frac{1}{|\L_{m_i}|} \sum_{x\in\L_{m_i}} \E\biggl(\int\nu^{\psi
_u}_{\L_{m_i}}[\xi](\ormd\phi)\bigl(\phi(e_1-x)-\phi(-x)\bigr)\biggr)\\
&&\qquad=\frac{1}{|\L_{m_i}|} \sum_{x\in\{\L_{m_i}\setminus\bar\Lambda
^{-m_i,m_i}_{-m_i}\}} \E\biggl(\int\nu^{\psi_u}_{\L_{m_i}}[\xi](\ormd
\phi
)\bigl(\phi(e_1-x)-\phi(-x)\bigr)\biggr)\\
&&\quad\qquad{}+\frac{1}{|\L_{m_i}|} \sum_{x\in\bar\Lambda^{-m_i,m_i}_{-m_i}}
\E
\biggl(\int\nu^{\psi_u}_{\L_{m_i}}[\xi](\ormd\phi)\bigl(\psi(e_1-x)-\phi
(-x)\bigr)\biggr).
\end{eqnarray*}
Most of the terms in the last equality in the above equation cancel and
we are left with
\begin{eqnarray*}
&&\E\biggl(\int\bar\mu_{m_i}^{\rho_u}[\xi](\ormd\eta)\eta(b_1)
\biggr)\\
&&\qquad=\frac{1}{|\L_{m_i}|}\biggl[ \sum_{x\in\bar\Lambda
^{-m_i,m_i}_{-m_i}}\psi(e_1-x)- \sum_{x\in\bar\Lambda
^{-m_i,m_i}_{m_i}}\bigl(u_1+\psi(-e_1-x)\bigr) \\
&&\hspace*{63pt}{}- \sum_{x\in\bar\Lambda^{-m_i,m_i}_{m_i}} \E\biggl(\int\nu
^{\psi_u}_{\L_{m_i}}[\xi](\ormd\phi)\bigl(\phi(-x)-\psi(-e_1-x)-u_1\bigr)
\biggr)\biggr]\\
&&\qquad=u_1+\frac{O(K,u_1)}{2m_i+1},
\end{eqnarray*}
uniformly in $m_i\in\N$, and where to bound the last term in the first
equality, we used Proposition~\ref{averaged}.
From this, it follows easily that we have, uniformly in $k\in\N$,
\[
\E\biggl(\int\hat\mu_k^{u}[\xi](\ormd\eta)\eta(b_1)\biggr)=u_1+\frac
{o(\log k)}{k}.
\]
%
Fix any large $M>0$. Then $ \eta(b) \wedge M \vee(-M)$ is bounded and
continuous, so for $\P$-a.s. all $\xi$, we have
\[
\lim_{k\rightarrow\infty}\int\hat\mu_k^u[\xi](\ormd\eta)\eta
(b) \wedge
M \vee(-M)=\int\hat\mu^u[\xi](\ormd\eta)\eta(b) \wedge M \vee(-M).
\]
Moreover, from Proposition~\ref{averaged} and Chebyshev's inequality,
we have
\[
\E\biggl(\int\hat\mu_k^u[\xi](\ormd\eta)\eta(b)\biggr)=\E\biggl(\int\hat
\mu_k^u[\xi](\ormd\eta)\eta(b) \wedge M \vee(-M)\biggr)+\frac{O(K)}{M^d},
\]
uniformly
in $k\in\N$. Therefore by sending $M$ to $\infty$, the convergence
of the
truncated $\eta$ together with the fact that $\int\hat\mu^u[\xi
](\ormd
\eta)\eta(b)$ is an integrable random variable, proves (\ref{tilt1}).
By symmetry, (\ref{tilt1}) holds for any $b\in(\zd)^*$.

To prove (\ref{intcond}), take any $b=(x_0,y_0)\in(\Z^d)^*$. Since
$(\phi(x_0)-\phi(y_0)-u\cdot(x_0 -y_0))^2\ge0$, by the weak
convergence of $({\hat{\mu}}_k^u)_{k\in\N}$ to ${\hat{\mu}}^u$
and by
Proposition~\ref{averaged}, we have
%
%
\begin{eqnarray}
\label{weakbound}
&&\E\biggl(\int{\hat{\mu}^u}[\xi](\mathrm{d}\eta)\bigl(\phi(x_0)-\phi(y_0)-u\cdot(x_0
-y_0)\bigr)^2\biggr)
\nonumber
\\[-8pt]
\\[-8pt]
\nonumber
&&\qquad\le\E\Bigl(\lim\inf_{k\rightarrow\infty}{\hat{\mu
}}_k^u\bigl(\phi(x_0)-\phi(y_0)-u\cdot(x_0 -y_0)\bigr)^2\Bigr)<K.
\end{eqnarray}
\upqed\end{pf*}

\begin{pf*}{Proof of Theorem \protect\ref{non-exist of Gibbs}} Suppose that
the infinite-volume gradient Gibbs measure does exist
and it satisfies
$\E|\int\mu[\xi](d\eta)V'(\eta(b))|<\infty\ \mbox{for all}\break
\mbox{bonds}\ b=(x,y)\in(\zd)^*$.
Then we have, in the present notation,
%
%
\begin{equation}
\label{eq:obenx}
\sum_{x\in\L}\xi(x)
= -\mathop{\sum_{x\sim y}}_{ x\in\L,y\in\L^c}
X_{(x,y)}[\xi]
\end{equation}
with $X_{b}[\xi]:=\int\mu[\xi](d\eta)V'(\eta(b))$ which was proved
in~\cite{EK}. We take $\L$ to be a~box,
divide both sides of the equation by $|\L|$ and take the limit $\L
\uparrow\Z^d$.
Then the right-hand side tends to zero if $d\geq1$,
while the left-hand side tends to the nonzero constant $\E(\xi(0))$ in
any dimension.
\end{pf*}
\subsubsection{Nonroughening in an averaged sense}


We will give next the following large deviation upper bound both for
the measures $\mu^{u}_{\L}[\xi]$, as defined in (\ref{1eqn1}), and for
the averaged measures ${\bar\mu}^{u}_{\L}[\xi]$, as defined in
(\ref{heldfixed}).
\begin{prop}\label{norafone} Suppose that $d\geq3$, $\E(\xi(0))=0$ and
$\E(\xi^2(0))<\infty$.

\begin{enumerate}
\item Then
there exist constants $K,\b, t_0>0$
such that for all but finitely many $N\in\N$, the following large
deviation upper bound
holds for all $t>t_0$ and for $\P$-almost all $\xi$:
%
%
\begin{eqnarray}
\label{norafone1}
&&\mu^{\rho_u}_{\L_N}[\xi]\biggl(
\frac{1}{2 |\L_N|}\sum_{x,y\in\L_N,|x-y|=1}\bigl(\phi(x)-\phi
(y)-u\cdot(x-y)\bigr)^2
>t \biggr)
\nonumber
\\[-8pt]
\\[-8pt]
\nonumber
&&\qquad \leq\exp(- \b|\L_N|t).
\end{eqnarray}

\item The same result holds for the averaged measures ${\bar\mu
}^{u}_{\L
}[\xi]$.
\end{enumerate}

\end{prop}

\begin{pf}
The assumption
$\E
(\xi^2(0))<\infty$
allows us to use the SLLN in Proposition~\ref{sllngree} along boxes
$\L
_N$ of side-length $N$, which implies that there exists
a nonrandom constant $K$ such that for $N$ large enough, we have $\frac
{1}{|\L_N| }\langle\xi, G_{\L_N} \xi\rangle_{\L_N}\leq K.$
Conditional on this bound,\vspace*{1pt}
one has by means of Lemma~\ref{Lemmaoneone} that $F_{\b,u,\L_N}[\xi
_{\L
_N}]\leq|\L_N| K$ (for a modified $K$)
which, by the exponential Chebychev inequality, implies the
concentration bounds of the form (\ref{norafone1}).\vadjust{\goodbreak}

To get \textit{the same type} of bounds for the measure ${\bar\mu
}^{u}_{\L
}[\xi]$, we need
to make use of the monotonicity in $\Lambda\in\Z^d$ of the quadratic
form $\langle\xi, G_{\L} \xi\rangle_{\L}$ stated in Proposition~\ref{monoform}.

Let us look at the quantity
\begin{eqnarray*}
&&\exp\hat F_{\b,u,\L}[\xi_{\L}]\\
&&\qquad:= \int\bar\nu^{u}_{\L}[\xi](
\ormd\phi)
\exp\biggl(
+\frac{\b}{2}\sum_{x,y\in\L_N,|x-y|=1}\bigl(\phi(x)-\phi(y)-u\cdot(x-y)\bigr)^2
\biggr)
\end{eqnarray*}
with the obvious definition for $\bar\nu^u_\L$. Note that we have the
following upper bound:
\[
e^{\hat F_{\b,u,\L}[\xi_{\L}]}\leq e^{\bar F_{\b,u,\L}} \frac
{1}{|\L
|}\sum_{x\in\L}e^{({\alpha }/{2})\langle\xi, G_{\L+x} \xi
\rangle_{\L+x}},
\]
by a straightforward application of the previous steps. By
Proposition~\ref{monoform} we have for each term
under the sum, the estimate $\langle\xi, G_{\L+x} \xi\rangle_{\L+x}
\leq\langle\xi, G_{\L+\L} \xi\rangle_{\L+\L}$
where $\L+\L:=\{x+y\dvtx  x,y \in\L\}$.
This gives us the estimate
\[
\label{1111}
\hat F_{\b,u,\L}[\xi_{\L}]\leq\bar F_{\b,u,\L} + \frac{\alpha
}{2}\langle
\xi, G_{\L+\L} \xi\rangle_{\L+\L}.
\]
From here the proof of the validity of the bounds stays the same.
\end{pf}

\section{Model B}\label{sec4}

The proof of Theorem~\ref{subadit} on surface tension for model B
follows the same argument as for model A, so it will be omitted. We
will focus instead on proving the existence of shift-covariant random
gradient Gibbs measures with given tilt. We consider the finite-volume
Gibbs measures
with tilt $u\in\R^d$ and boundary condition $\psi_u(x)=u\cdot x$ of
the form
\begin{eqnarray*}
&&\nu^{\psi_u}_{\L}[\omega](\ormd\phi)\\
&&\qquad=\frac{1}{Z_{\L}^{\psi
}[\omega
]}
\exp\biggl(
-\frac{1}{2}\mathop{\sum_{x,y\in\L} }_{|x-y|=1}V_{(x,y)}^\omega\bigl(\phi
(x)-\phi(y)\bigr)
\\
&&\hspace*{65pt}\qquad{}- \mathop{\sum_{x\in\L,y\in\partial\L}}_{|x-y|=1}V_{(x,y)}^\omega
\bigl(\phi
(x)-\psi(y)\bigr)
\biggr)\,\ormd\phi_{\L}\d_{\psi_u}(\ormd\phi_{\Z^d\setminus\L}).
\end{eqnarray*}
%
Similar to what we did for model A to prove tightness, we will consider
%
%
\begin{eqnarray}
&&\exp F_{\b,u,\L}[\omega_{\L}]
\nonumber\hspace*{-35pt}
\\[-2pt]
\\[-14pt]
\nonumber
&&\quad:= \int\nu^{\psi_u}_{\L}[\omega
](\ormd\phi)
\exp\biggl(
+\frac{\b}{2}\sum_{x,y\in\zd,|x-y|=1}\bigl(\phi(x)-\phi(y)-u\cdot(x-y)\bigr)^2
\biggr).\hspace*{-35pt}
\end{eqnarray}
By the same reasoning as for the proof of Lemma~\ref{Lemmaoneone}, we
get:
\begin{lem}
%
%
\begin{eqnarray}
F_{\b,u,\L}[\omega_{\L}]&\leq&
-|\L|\bigl(\s_{\L}^{A-\b}[\omega=0](u=0) - \s_{\L}^{C_2}[\omega=0](u=0)
\bigr)\nonumber\\
&&{}+\mathop{\sum_{x,y\in\L\cup\partial\L}}_{|x-y|=1}B_{(x,y)}^\omega
\nonumber
\\[-8pt]
\\[-8pt]
\nonumber
&&{}-\frac{A-\b-C_2}{2} \mathop{\sum_{x,y\in\L\cup\partial\L}}_{|x-y|=1}\bigl((x-y)\cdot u\bigr)^2\\
&=:&\tilde F_{\b,u,\L}+ \mathop{\sum_{x,y\in\L\cup
\partial\L}}_{ \atop|x-y|=1}B_{(x,y) }^\omega,\nonumber
\end{eqnarray}
where the first term on the right-hand side is a nonrandom quantity
which is bounded by a constant
times $|\Lambda|$.
\end{lem}

Note that the critical dimension for existence changes from $d=3$, as
it was in model A, to $d=1$. The reason for this change is the absence
of the term $\langle\xi, G_{\L} \xi\rangle_{\L}$ in the formula for
$F_{\b,u,\L}[\omega_{\L}]$ above, and which term, present in the
formula for $F_{\b,u,\L}[\xi_{\L}]$ in model A, diverges for $d=2$ when
averaged over the disorder.

Define $\mu^{\rho_u}_\L[\omega]$ and $\bar\mu^{\rho_u}_\L
[\omega]$ as
for model A.
As in Proposition~\ref{averaged} from model~A, we have the following
result on the tightness of the family of finite-volume random $\nabla
\phi$-Gibbs measures $\mu^{u}_{\L_N}[\omega]$ averaged over the disorder.
\begin{prop}
\label{averaged1}
Suppose that $d\geq1$. Then there exists a constant \mbox{$K>0$} such that
for all bonds
$x_0, y_0\in\Z^d$, with $|x_0-y_0|=1$, we have that the measure
$P^u_{\L}(\mathrm{d}\phi):=\int\P(\mathrm{d}\omega)\mu^{\rho_u}_{\L}[\omega
](\mathrm{d}\phi)$
satisfies the estimate
%
%
\begin{equation}
\limsup_{N\uparrow\infty} P^u_{\L_N}\bigl(\phi(x_0)-\phi(y_0)\bigr)^2 \leq K.
\end{equation}
Hence the sequence of measures $P^u_{\L_N}$ is tight and thus possesses
a disorder-independent
limit measure (along subsequences of volumes) on gradient configurations.
\end{prop}

\begin{pf} We proceed exactly as for model A to get the bound
%
%
\begin{eqnarray}
P^{u}_{\L}(f)
&\leq&\frac{2}{\b|\L|} \E\log\mu^{\rho_u}_{\L}[\omega]
\nonumber
\\[-8pt]
\\[-8pt]
\nonumber
&&{}\times\biggl(\exp
\biggl(\frac{\b}{2}\sum_{x,y\in\zd,|x-y|=1}\bigl(\phi(x)-\phi(y)-u\cdot
(x-y)\bigr)^2\biggr)\biggr),
\end{eqnarray}
which gives us
%
%
\begin{equation}
P^{u}_{\L}(f)
\leq\frac{2}{\b|\L|}\tilde F_{\b,u,\L}+\frac{2}{\b|\L|}\biggl(\mathop{\sum_{x,y\in\L\cup\partial\L}}_{|x-y|=1}\E B^\omega_{(x,y) }
\biggr),
\end{equation}
which is bounded uniformly in $\L$.
\end{pf}


Theorem~\ref{thm1}(b) follows now immediately from Proposition
\ref{averaged1} by similar reasoning as in the proof of Theorem~\ref{thm1}(a).

Similar to the proof of Proposition~\ref{norafone}, we have the
following large deviation upper bound for the finine volume Gibbs
measures $\mu^{\rho_u}_{\L}[\omega]$ and $\bar\mu^{u}_{\L
}[\omega]$.

\begin{prop} Suppose that $d\geq1$. Then
there exist constants $K,\b,\allowbreak t_0>0$ such that for all realizations
$\omega\in\Omega$ and
for all $N\in\N$ the following large deviation upper bound
holds for all $t>t_0$:
%
%
\begin{eqnarray}
&&\mu^{\rho_u}_{\L_N}[\omega]\biggl(
\frac{1}{2 |\L_N|}\sum_{x,y\in\L_N,|x-y|=1}\bigl(\phi(x)-\phi
(y)-u\cdot(x-y)\bigr)^2
>t \biggr)
\nonumber
\\[-8pt]
\\[-8pt]
\nonumber
&&\qquad\leq\exp(-\b|\L_N|t)
\end{eqnarray}
and
%
%
\begin{eqnarray}
&&\bar\mu^{u}_{\L}[\omega]\biggl(
\frac{1}{2 |\L_N|}\sum_{x,y\in\L_N,|x-y|=1}\bigl(\phi(x)-\phi
(y)-u\cdot(x-y)\bigr)^2
>t \biggr)
\nonumber
\\[-8pt]
\\[-8pt]
\nonumber
&&\qquad \leq\exp(- \b|\L_N|t).
\end{eqnarray}
\end{prop}


\begin{appendix}
\section*{Appendix}\label{app}

\subsection{Why the Gibbs measure does not exist for model A in $d=3,4$
for $V(s)=s^2/2$}

We will prove next that for model A in $d=3,4$, there\break exists no
infinite-volume Gaussian Gibbs measure with $s:=\break\sup_{x\in\zd}\E
|\int\nu[\xi](\ormd\phi)\phi(x)|<\infty$. Take $\L
_N:=[-N,N]^d\cap
\Z^d, N\in\N,$ and let $\psi\in\R^{\zd}$ be an arbitrary boundary
condition. Then we have for the finite-volume Gibbs measure
%
%
\setcounter{equation}{0}
\begin{equation}
\label{paininthearse}
\int\nu_{\L_N}^\psi[\xi](\ormd\phi)\phi(0)=\sum_{z\in{\L
_N}}G_{\L
_N}(0,z)\xi(z)+
{\E}_0 (\psi(X_{\t_{{\L_N}}})).
\end{equation}
Here the expectation ${\E}_0$ is w.r.t. a nearest-neighbor random walk
$X:=(X_k)_{k\in\N}$ started at $0$ with Green's function $(G_{\L
_N}(0,y))_{y\in{\L}_N}$,
and the second term is what we obtain for the nondisordered model. We
defined $\tau_{\L_N}:=\inf\{k\ge0\dvtx X_k\in{\L_N}^c\}$, so $X_{\tau
_{\L
_N}}$ is the position of the random walk when it exits $\L_N$.
Suppose that there is a random infinite-volume Gibbs measure~$\nu[\xi]$
in $d=3,4$.
Average (\ref{paininthearse})\vadjust{\goodbreak} over the boundary conditions $\psi$
w.r.t. the measure~$\nu[\xi]$ and use the DLR equation to conclude that
%
%
\begin{equation}\qquad
\label{194}
\int\nu[\xi](\ormd\phi)\phi(0)=\sum_{z\in{\L_N}}G_{\L
_N}(0,z)\xi(z)+
{\E}_0 \int\nu[\xi](\ormd\phi)(\phi(X_{\t_{\L_N}})).
\end{equation}
The expectation under the disorder for the second term in (\ref{194})
stays bounded uniformly in $\Lambda_N$ under our hypothesis; in fact,
we have
%
%
\begin{eqnarray}
\label{204}
&&\E\biggl|{\E}_0 \int\nu[\xi](\ormd\phi)(\phi(X_{\t_{{\L}_N}}))\biggr|\nonumber\\[-2pt]
&&\qquad=
\E\biggl|\sum_{u \in\partial{\L}_N}{\P}_0(X_{\t_{{\L}_N}}=u)\int\nu
[\xi
](\ormd\phi)\phi(u)\biggr|\\[-2pt]
&&\qquad\leq\sum_{u \in\partial{\L}_N}{\P}_0(X_{\t_{{\L}_N}}=u)\E
\biggl|\int
\nu[\xi](\ormd\phi)\phi(u)\biggr|\leq s.\nonumber
\end{eqnarray}
The left-hand side of (\ref{194}) is a proper random variable and $
({\E}_0 \int\nu[\xi](\ormd\phi)\times (\phi(X_{\t_{{\L}_N}})))_{{\L
}_N\subset\Z^d}$
is a tight family of random variables by (\ref{204}). However, $
(\sum_{z\in{\L}_N}G_{{\L}_N}(0,z)\xi(z))_{{\L}_N\subset\Z^d}$ is
not a tight family because a simple characteristic function calculation
shows that
\[
\frac{\sum_{z\in{\L}_N}G_{{\L}_N}(0,z)\xi(z)}{\sqrt{\sum_{z\in
{\L
}_N}G^2_{{\L}_N}(0,z)}}
\]
converges to a standard normal as $N\uparrow\infty$, since $\sum
_{z\in
{\L}_N}G^2_{{\L}_N}(0,z)$ diverges in $d=3,4$. This leads to a
contradiction in (\ref{194}) as ${\L}_N\uparrow\Z^d$.

The identity in (\ref{paininthearse}) is based on exact computations
for multivariate Gaussian distributions, which we do not have for
nonquadratic potentials.
For the more general class of potentials satisfying (\ref{tag2}) and
(\ref{tag22}), we expect the conclusion to be the same.

\subsection{Why the Brascamp--Lieb inequality does not solve the problem}
A~different route to proving the existence of random gradient Gibbs
measures uses the Brascamp--Lieb inequality.
It states that for $\g$ a centered Gaussian distribution on $\R^d$
and a distribution $\mu$ on $\R^d$ such that there exists $d\mu
/d\gamma
= e^{-f}$ for a convex function $f$,
one has for all $v\in\RR^d$ and for all convex real functions
$F$, bounded below, that
%
\begin{equation}
\label{111}
\mu\bigl(F\bigl(v \cdot\bigl(X -\mu(X)\bigr)\bigr)\bigr)\leq\g\bigl(F(v \cdot
X)\bigr).
\end{equation}
The above is the formulation by Funaki in~\cite{F2}. An application of
(\ref{111}) to our disorderd case would give, for example, that
%
%
\begin{eqnarray}
&&\mu^{\rho_u}_{\L}[\x] \bigl( \bigl[
\phi(x_0)-\phi(y_0)- \mu^{\rho_u}_{\L}[\x] \bigl( \phi(x_0)-\phi
(y_0)\bigr)\bigr]^2\bigr)
\nonumber
\\[-9pt]
\\[-9pt]
\nonumber
&&\qquad\leq\g_{\L}\bigl( [\phi(x_0)-\phi(y_0)]^2\bigr),
\end{eqnarray}
where $\g_{\L}$ is the corresponding Gaussian measure.
The right-hand side is uniformly bounded\vadjust{\goodbreak} in $\L$, so that would prove
a.s. tightness for strictly convex potentials $V$
if we can prove that the expected values of the local tilts of the
interface taken over the Gibbs distribution have
limits for almost surely every realization of disorder, that is, if we
can prove that
%
%
\begin{equation}
\label{122}
\lim_{|\L|\uparrow\infty}\mu^{\rho_u}_{\L}[\x] \bigl( \phi
(x_0)-\phi
(y_0)\bigr)
\end{equation}
exists a.s. for $x_0,y_0in\L,$ with $|x-y|=1$. However, currently we
do not have~a~way either to prove (\ref{122}) or to prove the existence
of the\break $\lim_{|\L|\uparrow\infty}{\bar{\mu}}^{\rho_u}_{\L}[\x] (
\phi(x_0)-\phi(y_0))$, as introduced in (\ref{heldfixed}), in the
presence of disorder. Note that in the model without disorder, we can
show for strictly convex potentials $V$ the existence of the last limit
by Brascamp--Lieb inequality coupled with shift-invariance arguments.
\end{appendix}

\section*{Acknowledgments}
We thank David Brydges for pointing out to us a~reference for
Proposition~\ref{monoform},
and Marek Biskup, Jean-Dominique Deuschel and Marco Formentin for
stimulating discussions. We also thank Noemi Kurt, Rongfeng Sun and two
anonymous referees for very useful comments, which greatly improved the
presentation of the manuscript.

%


\printaddresses

\end{document}